\documentclass[twocolumn]{autart}    

\usepackage{graphicx}          
\usepackage{cite}
\usepackage{amsmath,amssymb,amsfonts}
\usepackage{algorithm}
\usepackage{algorithmicx}
\usepackage{algpseudocode}
\usepackage{graphicx}
\usepackage{hyperref}
\hypersetup{hidelinks=true}
\usepackage{textcomp}

\usepackage{enumerate}
\usepackage{caption,subcaption}
\usepackage{float,subfloat,stfloats}
\usepackage{array}
\usepackage{booktabs} 
\usepackage{threeparttable}

\begin{document}

\begin{frontmatter}

\title{Distributed Identification of Stable Large-Scale Isomorphic Nonlinear Networks Using Partial Observations\thanksref{footnoteinfo}} 

%
\thanks[footnoteinfo]{This work was supported by the National Natural Science Foundation of China (Grant No. 61991414, 62088101, 6193000461), Chongqing Natural Science Foundation CSTB2023NSCQ-JQX0018, and Beijing Natural Science Foundation L221005. The material in this paper was not presented at any conference.}
\thanks[*]{Corresponding author Chengpu Yu.}
\author[Beijing]{Chunhui Li}\ead{chunhui\_li@bit.edu.cn},    
\author[Beijing]{Chengpu Yu\thanksref{*}}\ead{yuchengpu@bit.edu.cn}     
%
\address[Beijing]{School of Automation, Beijing Institute of Technology, Beijing, 100081, China}  


\begin{abstract}                          
Distributed parameter identification for large-scale multi-agent networks encounters challenges due to nonlinear dynamics and partial observations. Simultaneously, ensuring the stability is crucial for the robust identification of dynamic networks, especially under data and model uncertainties. To handle these challenges, this paper proposes a particle consensus-based  expectation maximization (EM) algorithm. The E-step proposes a distributed particle filtering approach, using local observations from agents to yield global consensus state estimates. The M-step constructs a likelihood function with an a priori contraction-stabilization constraint for the parameter estimation of isomorphic agents. Performance analysis and simulation results of the proposed method confirm its effectiveness in identifying parameters for stable nonlinear networks.
\end{abstract}

\end{frontmatter}

\section{Introduction}
\par Large-scale networks play a crucial role in contemporary engineering areas, with applications ranging from intelligent traffic management to environmental monitoring and drone swarm coordination. These applications generate, transmit, and utilize vast data sets, yet the limited communication resources prevent agents from accessing the global information of the network. Moreover, the dynamics guiding agents in task execution are variable, intricate, and challenging to accurately characterize. These complexities hinder precise task execution, emphasizing the need for accurate identification of unknown dynamics using local inputs and outputs.

\par With the rapid expansion of networks, centralized parameter identification methods face challenges related to computational efficiency and data storage \cite{6758360,YU2019108517,HUANG2022110261}. The presence of unmodeled dynamic interactions among agents adds complexity to centralized approaches\cite{article,ZHAO2024108735}, attracting attention to distributed parameter identification methods. Various approaches have been explored, such as the distributed solution proposed by Chan et al. \cite{8675360}, which employed the augmented Lagrangian for nonlinear optimization problems in networked systems with partially separable functions and nonlinear constraints. This approach achieved robust estimation, particularly for problems with equality constraints. Revay et al. \cite{9547745} achieved distributed identification of large-scale networks, incorporating convex constraints on stability and monotonicity using the alternating direction method of multipliers. Gao et al. \cite{Gao} addressed the challenging problem of autonomously inferring complex network dynamics from empirical data, presenting a robust two-stage method successfully validated across various real and synthetic networks. Despite these advancements, current distributed dynamics identification methods for large-scale networks, particularly those based on nonlinear models, still grapple with challenges such as the need for extensive learning data and ideal conditions. Addressing issues like non-Gaussian noise, incomplete observations, and model uncertainties require further research.

\par In practical networks, agents often exhibit complex nonlinear dynamics that demand a flexible framework to capture the intricate interactions and dependencies among components \cite{ZHANG202325,MASTI2021109666,GEDON2021481,HUANG2024111349}. The observed data from these networks typically includes only inputs and noisy outputs, while internal states and external disturbances remain unobservable. To address the identification of a single nonlinear state-space model with skewed measurement noise, Liu et al.\cite{9843908} proposed an EM algorithm for simultaneous estimation of unknown state and model parameters. In addition, Ramadan et al. \cite{RAMADAN2022110482} iteratively combined particle filtering (PF) and the EM algorithm, deriving a maximum likelihood recursive state estimator for nonlinear state-space models. While traditional EM algorithms have successfully addressed parameter identification problems in systems with hidden variables, their distributed application in large-scale networks with internal state interactions remains challenging. Current methods have not effectively tackled issues like incomplete state measurements and unknown interacting functions in large-scale networks, indicating the need for further exploration and advancement.

\par A significant challenge in nonlinear modeling is the absence of stability guarantees, where models lacking stability may fit training data well but exhibit unpredictable behavior for unseen inputs \cite{8401713}. Ensuring stability becomes crucial, particularly when the structural details of the system are not fully known. Scholars have explored methods to ensure a priori stability for nonlinear models. Manchester et al.\cite{MANCHESTER2012328} presented a nonlinear system identification method that minimizes the convex upper bound of the long-term simulation errors within a convex set of stable nonlinear models. Tobenkin et al.\cite{7907229} addressed the simulation error minimization, which is non-convex with respect to the model parameters, using Lagrangian relaxation, dissipation inequalities, contraction theory, and semidefinite programming. Umenberger et al.\cite{8447503} developed a path-following interior-point algorithm, leveraging special problem structures to reduce computational complexity from cubic to linear growth with the dataset length. The above-mentioned algorithms mainly focus on cost functions and optimization algorithms, but there remains an unexplored frontier in leveraging a priori nonlinear stable models to address unmodeled dynamics in large-scale networks and limited data for learning system dynamics. Although recent results \cite{9547745} provided stability guarantees for models of networks, ongoing research on methods guaranteeing a priori stability for nonlinear models in large-scale networks deserves attention.

\par This study addresses the challenges of distributed identification in large-scale networks characterized by nonlinear dynamics and partial state observations by developing a particle consensus-based distributed particle
EM (PC-DPEM) algorithm. Additionally, it introduces contraction stability constraints to ensure the robustness of nonlinear dynamics. The main contributions of this paper can be summarized as follows:

\begin{enumerate}
	
	\item {The particle consensus-based distributed particle filtering (PC-DPF) approach is proposed under the EM framework to extend the identification of the single system with latent variables in previous research \cite{SCHON201139} to large-scale networked systems with nonlinear dynamic interactions. It addresses the challenge of distributed identification of isomorphic large-scale nonlinear networks using only local and partial state observations based on sparse topological communication.}
	
	\item {For partially observable agent dynamics in a large-scale nonlinear network, a priori contraction stability constraint is introduced to ensure the robust parameter identification in the presence of unknown measurement noise and uncertainties in the interaction functions among neighboring agents. Compared with the recent study in \cite{Gao}, the stability-constrained identification method in this paper can yield reliable results using partial state observations with measurement noises.}
	
\end{enumerate}

\par Subsequent sections provide a detailed exposition of relevant preliminaries and the problem setup (Section \ref{Sec2}), the EM framework for the identification of large-scale nonlinear networks (Section \ref{Sec3}), distributed state estimation by proposing a PC-DPF  algorithm (Section \ref{Sec4}), parameter estimation for the state-space represented isomorphic network (Section \ref{Sec5}), performance analysis of the proposed PC-DPEM algorithm (Section \ref{Sec6}), and simulation validation of the proposed identification algorithm in identifying complex large-scale networks (Section \ref{Sec7}). Finally, conclusions are drawn in Section \ref{Sec8}.

\section{Preliminaries and Problem Setup}
\label{Sec2}

\subsection{Network Topology and Communication}

\par The considered multi-agent nonlinear network consists of $V$ agents with limited computational and communication capabilities, and this directed network has $E$ edges that allow self-looping connections. The directed network topology is represented by $\mathcal{G}=\{\mathcal{V},\mathcal{E}\}$, where $\mathcal{V}=\{1,\ldots,V\}$ is the set of agents, and $\mathcal{E}\subseteq \mathcal{V}\times \mathcal{V}$ denotes the set of directed edges satisfying $|\mathcal{E}|=E$. Each agent $v$ has its own state $x_v$, input signal $u_v$, and output measurement $y_{v}$, where $v\in\mathcal{V}$. These data can be transmitted along the directed edge $e_{vj}\in\mathcal{E}$ if there is a directed link from agent $v$ to agent $j$, enabling adjacent agents to achieve local data sharing.

\par The topological relationships among different agents in a multi-agent network can be compactly expressed in a $V \times V$ edge-agent adjacency matrix $\mathcal{A}$\cite{6425904}. Each row represents an edge originating from a particular agent $v$, and each column represents the agent to which the edge ends. This edge-agent adjacency matrix also describes the self-looping connections in the network topology. The neighbors of agent $v$ defined by the edge-agent adjacency matrix $\mathcal{A}$ are categorized into two sets: predecessor neighbors and successor neighbors, denoted as $\mathcal P_v=\{j|e_{jv}\in \mathcal{E}, j\leqslant v\}$ and $\mathcal{S}_v=\{j|e_{vj}\in \mathcal{E}, v < j\}$, respectively. For an edge $e_{vj}$, the value of $\mathcal{A}_{e_{vj}}$ can be defined as
\begin{equation}
	\label{eq3_A}
	\setlength{\abovedisplayskip}{3pt}
	\setlength{\belowdisplayskip}{3pt}
	\textstyle
	\begin{aligned}
		\mathcal{A}_{e_{vj}}=\left\{
		\begin{array}{ll}
			1\quad &\text{if}\; j\in\mathcal{S}_v,\\
			0\quad &\text{otherwise}.
		\end{array}
		\right.
	\end{aligned}
\end{equation} 
It should be noted that the neighbors of agent $v$ are free to join or leave.

\par Utilizing this directed network topology, each agent engages in iterative information exchange through a direct communication mechanism. Agents share and integrate data with their neighbors, collaboratively adjusting parameter estimates based on local observations, which enables to enhance the accuracy and consistency of dynamic parameter estimation of the isomorphic network.

\subsection{Network Nonlinear State Space Model}

\par Large-scale networks involve dynamic interactions among individual agents, and unknown interaction functions will pose challenges for accurate network identification. For an isomorphic networked system with the described communication topology, this study constructs a quadruple model class $\mathcal{M}=(f,h,g,\theta)$ for each agent $v$. In this model class, $f : \mathbb{R}_n \times \mathbb{R}_m \times \mathbb{R}_s \times \mathbb{R}_{q} \to \mathbb{R}_n$, $h : \mathbb{R}_n \times \mathbb{R}_m \times \mathbb{R}_w \times \mathbb{R}_{q} \to \mathbb{R}_p$, and $g : \mathbb{R}_n \times \mathbb{R}_n \to \mathbb{R}_n$ are continuously differentiable functions, while $\theta\in \mathbb{R}_q$ represents a q-dimensional vector, indicating unknown parameters embedded in $f$ and $h$. 

\par The considered network model is represented by a discrete-time nonlinear state-space model as 
\begin{subequations}
	\label{eq1:combined}
	\setlength{\abovedisplayskip}{3pt}
	\setlength{\belowdisplayskip}{3pt}
	\begin{align}
		\setlength{\abovedisplayskip}{3pt}
		\setlength{\belowdisplayskip}{3pt}
		\textstyle
		{x}_v(t+1)=&\textstyle f(x_v(t),u_v(t),\varepsilon_v(t),\theta)\label{1a}\\
		&\textstyle+\sum_{j\in\mathcal P_v}g_j(x_j(t)-x_v(t)), \notag\\
		y_v(t)=&\textstyle h(x_v(t),u_v(t),\eta_{v}(t),\theta), v\in\mathcal{V}.\label{1b}
	\end{align}
\end{subequations}				
The transition equation (\ref{1a}) captures the state evolution of agent $v$ denoted by $x_v \in \mathbb{R}_n$ over time. Here, $f$ is a nonlinear mapping function, $u_v\in \mathbb{R}_m$ is the exogenous input, and $\varepsilon_v \in \mathbb{R}_s$ represents independent process noise. The current state $x_v(t)$ of agent $v$ may be influenced by its neighbors $j$, where $j\in\mathcal P_v$. These influences, represented by $g_j$, are defined as the interaction function of state differences against the neighboring agents. By linearizing around $x_v(0)$, assume that $\sum_{j\in\mathcal P_v}g_j(x_j(t)-x_v(t))$ satisfies $\sum_{j\in\mathcal P_v}G(x_v)=\sum_{j\in\mathcal P_v}\frac{\partial}{\partial x_v}{{g}_{j}}(x_j(t)-x_v(t))\leqslant\kappa$ for some $\kappa>0$. This assumption is motivated by the conditions established in Definition \ref{def1} and Theorem \ref{thm4}. The measurement equation (\ref{1b}) represents the relationship between the state $x_v$ and measurement $y_{v} \in \mathbb{R}_p$. Here $h$ is a nonlinear mapping function, and $\eta_{v} \in \mathbb{R}_w$ represents independent measurement noise. When $p<n$, it is called partial state observation, which will be investigated in this paper.

\par Considering the random noises $\varepsilon_v$ and $\eta_v$, the model \eqref{eq1:combined} can be represented in stochastic form as
\begin{subequations}
	\label{eq2:combined}
	\setlength{\abovedisplayskip}{3pt}
	\setlength{\belowdisplayskip}{3pt}
	\begin{align}
		\setlength{\abovedisplayskip}{3pt}
		\setlength{\belowdisplayskip}{3pt}
		\textstyle
		x_v(t+1)&\sim p_{\theta}(x_v(t+1)|x_{\mathcal P_v}(t)), \label{2a} \\ 	 
		y_{v}(t)&\sim p_{\theta}(y_{v}(t)|x_{\mathcal P_v}(t)),  \label{2b}
	\end{align}
\end{subequations}
where $v\in\mathcal{V}$. Here, (\ref{2a}) represents the probability density function (pdf) of the state $x_v(t+1)$ given $x_{\mathcal P_v}(t)$, $u_v(t)$, and $\theta$. Similarly, (\ref{2b}) represents the pdf of the measurement $y_{v}(t)$. For brevity, the input signal $u_v(t)$ is omitted when no confusion will be raised.

\subsection{Problem of Interest}

\par This study aims to identify the unknown model parameter $\theta$ in the dynamics (\ref{eq1:combined}) for each agent $v$, $v\in\mathcal{V}$ within an isomorphic network, where the parameter $\theta$ are identical across different agents. To deal with this problem, the network topology represented by the adjacency matrix $\mathcal{A}$ is perfectly known. The nonlinear mapping functions $f$ and $h$ in (\ref{eq1:combined}) are linearly parametrized by $\theta$ as $f(\cdot)=\sum_{i=0}^q \theta_i f_i(\cdot)$ and $h(\cdot)=\sum_{i=0}^q \theta_i h_i(\cdot)$, where the basis functions $f_i$ and $h_i$ are typically chosen as multivariate polynomials, trigonometric polynomials, etc. Note that the interaction function $g$ is unavailable, which poses challenges to the accurate parameter identification. Additionally, $\varepsilon_v$ and $\eta_v$ represent independently and identically distributed (i.i.d) noises following a known form (e.g., Gaussian), but the parameters of the noise distribution, such as mean and variance, are unknown and can be absorbed into $\theta$ for identification purposes.

\begin{figure}   
	\centering
	\includegraphics[width=0.45\textwidth,trim=0 0 0 0,clip]{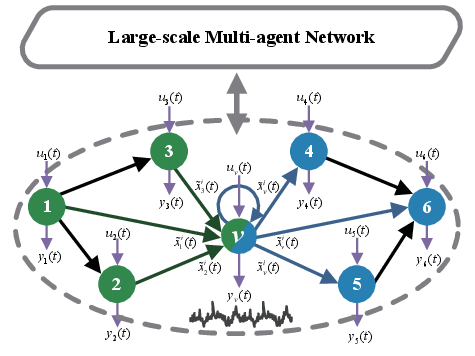}
	\caption{Local communication topology of a large-scale network. $\{\tilde x_j^i(t)\}_{j=1}^3$ are the local state estimation particles transmitted to agent $v$ from its predecessors, $\{\tilde{x}_v^i(t)\}$ are the local state estimation particles transmitted from agent $v$ to its successors. $\{\{u_i(t)\}_{i=1}^6, u_v(t)\}$ and $\{\{y_i(t)\}_{i=1}^6, y_v(t)\}$ are locally measurable inputs and outputs of agent $v$ respectively.}
	\label{f111}
\end{figure}

\par The parameter identification relies on a batch of $T$ input and output measurements for each agent, which are denoted as $u_v = [u_v(1), \ldots, u_v(T)]$ and $y_{v} = [y_{v}(1), \ldots, y_{v}(T)]$ for $v\in\mathcal{V}$. In subsequent discussions, the input and output measurements of agent $v$ are represented respectively as $U_v(1:T)$ and $Y_{v}(1:T)$, which are abbreviated as $U_v(T)$ and $Y_{v}(T)$. Fig. \ref{f111} illustrates a local network consisting of agent $v$ and its neighbors $\mathcal P_v\bigcup\mathcal S_v$, where the local network topology is known and the input-output responses are measurable in the presence of noise.

\begin{figure*}   
	\centering
	\includegraphics[width=\textwidth,trim=0 0 0 0,clip]{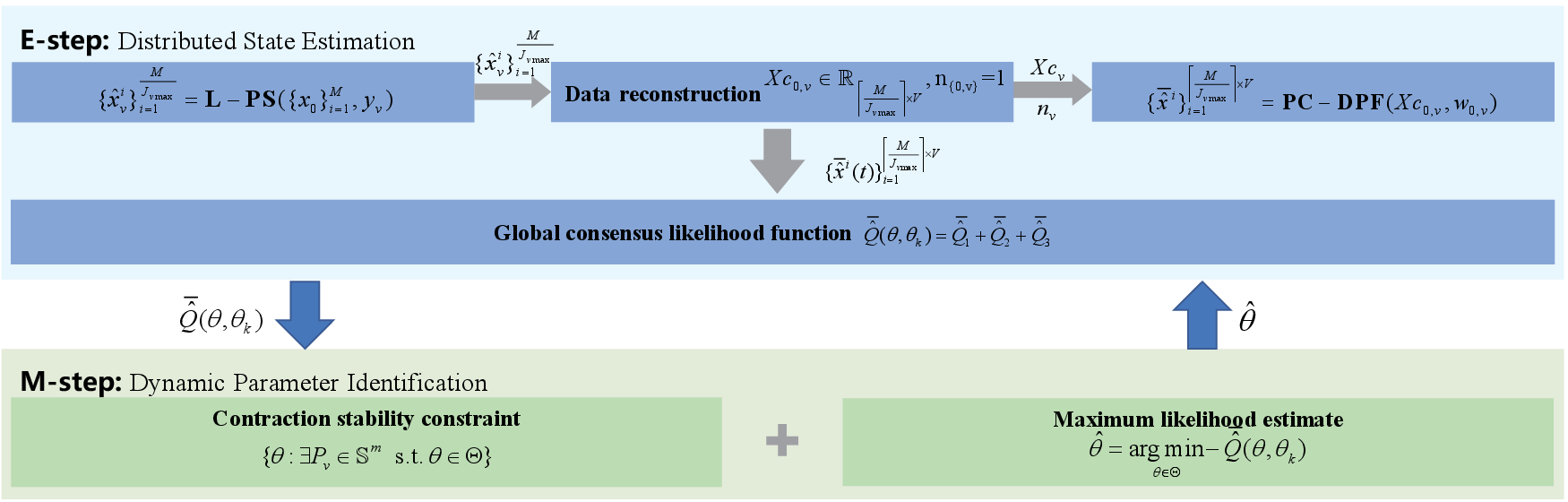}
	\caption{The implementation strategy for the PC-DPEM algorithm, where L-PS is a shorthand form of individual agent particle smoothing and PC-PDF is a shorthand form of particle consensus-based distributed particle filtering.}
	\label{ff}
\end{figure*}

\par The main challenges of this research involve two aspects. Firstly, there is a need to address the issue of incomplete local information, i.e., agents can only obtain local measurements, while the state is not fully measured. The fusion of network information becomes crucial to achieve consistent and accurate parameter identification results for isomorphic networks. In addition, the bilinear coupling phenomenon between unknown states and network parameters is difficult to handle for the concerned identification problem with latent variables. Secondly, it is necessary to tackle the challenges arising from unknown noise statistics and unknown dynamic interactions among agents. These two unknown terms may lead to system instability or unreliable parameter estimation. To ensure a predictable response of each agent under uncertain conditions, a contraction stability concept is appropriately applied\cite{7907229}, defined as
\begin{defn}
	\label{def1}
	(Global Incremental $l_2$ Contraction Stability): The deterministic part of the transition equation (\ref{1a}) in model (\ref{eq1:combined}), defined as ${x}_v(t+1)=f(x_v(t),u_v(t),\varepsilon_v(t),\theta)$, is considered incrementally $l_2$ contraction stable if, for any two distinct initial values ${x}_{v1}(0)$, ${x}_{v2}(0)$ and every input sequence $u_v(t)$, $v\in\mathcal{V}$, $t=1,2,\ldots$, the states ${x}_{v1}(t)$ and ${x}_{v2}(t)$ satisfy $\sum_{t=0}^{\infty}|{{x}_{v1}}(t)-{{x}_{v2}}(t){{|}^{2}}<\infty$.
\end{defn}
Stability can be proved by finding a differential storage function $V: \mathbb{R}^{n} \times \mathbb{R}^{n} \to \mathbb{R}$. The function $V$ can be any positive semidefinite function with respect to $x_v(t)$ and its differential $\Delta_v(t)$ in the initial conditions, as well as some fixed $\kappa>0$, there is a difference dissipation inequality satisfying  $V(x_v(t+1), \Delta_v(t+1)) - V(x_v(t), \Delta_v(t)) \leqslant - \kappa|\Delta_v(t)|^2$. 

\par To deal with the network identification problem, this study proposes a PC-DPEM algorithm within the framework of EM, as illustrated in Fig. \ref{ff}, to accomplish distributed parameter identification. The algorithm operates iteratively in two main steps:

\begin{enumerate}
\item{In the E-step, the PC-DPF algorithm is employed for global consensus state estimation. Subsequently, a distributed particle approximation function $\bar{\hat{Q}}({\theta}, \theta_{k})$ is computed locally for each agent, which serves as an approximation of the global likelihood function.}
\item{In the M-step, optimization of $\bar{\hat{Q}}({\theta}, \theta_{k})$ takes place under the contraction stability constraint, which is crucial for achieving distributed parameter estimation.}
\end{enumerate}

\section{The EM Algorithm for Networks}
\label{Sec3}

\par Since each agent in an isomorphic multi-agent network can only acquire local measurements and the agent's internal states are not fully observable, this study considers the application of an EM framework for network dynamic parameter identification. To facilitate the later description, (\ref{eq1:combined}) can be compactly represented as
\begin{equation}
	\label{eq3.1.1}
	\setlength{\abovedisplayskip}{3pt}
	\setlength{\belowdisplayskip}{3pt}
	\begin{aligned}
		\setlength{\abovedisplayskip}{3pt}
		\setlength{\belowdisplayskip}{3pt}
		\textstyle
		{x}(t+1)&\textstyle=f(x(t),\varepsilon(t),\theta )+{g}(x(t)-\mathcal{A}x(t)), \\
		y(t)&\textstyle=h(x(t),\eta(t),\theta),
	\end{aligned}
\end{equation}
where $x$, $\varepsilon$, $\theta$, $y$, and $\eta$ are vectors composed of $x_v$, $\varepsilon_v$, $\theta_v$, $y_v$, and $\eta_v$ for $v \in\mathcal{V}$, respectively. The structures of the differentiable nonlinear mapping functions $f(\cdot)$, $h(\cdot)$, and $g(\cdot)$ are consistent with the distributed representation in (\ref{eq1:combined}).

\par The EM algorithm treats the unobserved states $X(T) = [x(1), \ldots, x(T)]$ as latent variables and constructs a joint likelihood function which is denoted as $L_{\theta}(X(T), Y(T)) = \log p_{\theta}(X(T), Y(T))$. By applying the principle of maximization-minimization \cite{UMENBERGER2018280}, the algorithm iteratively seeks the maximum likelihood parameter estimate $\hat{\theta}$.

\par Each iteration of the EM algorithm involves two steps. In the E-step, due to the existence of the Theorem \ref{thm1}, the algorithm computes the minimum variance estimate $Q({\theta}, \theta_{k})$ of $L_{\theta}(X(T),Y(T))$, denoted as 
\begin{equation}
	\label{eq3.1.2}
	\setlength{\abovedisplayskip}{3pt}
	\setlength{\belowdisplayskip}{3pt}
	\textstyle
	\begin{aligned}
		Q({\theta}, \theta_{k}):=&\textstyle E_{{\theta}_{k}}[L_{\theta}(X(T), Y(T)) | Y(T)] \\
		=&\textstyle \int L_{\theta}(X(T), Y(T))\cdot \\	
		&\textstyle \quad p_{\theta_{k}}(X(T)|Y(T)) dX(T),
	\end{aligned}
\end{equation}
to tackle the unknown latent variable $X(T)$, where $E_{\theta_k}$ denotes the expectation operator with respect to the distribution $p_{\theta_k}(\cdot)$. Subsequently, in the M-step, this function is maximized to obtain a new estimate of the parameter vector $\theta$. The EM algorithm iteratively updates the parameter estimate $\hat{\theta}_{k}$ so as to maximize the log-likelihood function $L_{\theta}(Y(T))$\cite{GIBSON20051667}.

\begin{thm}
	\label{thm1}
	 Suppose that $\theta_{k+1}=\arg\max_{\theta} Q(\theta,\theta_k)$. Then, $\theta_{k+1}$ will result in an increase of the log-likelihood function, i.e., $L_{\theta_{k+1}}(Y(T)) > L_{{\theta}_{k}}(Y(T))$.
\end{thm}
\begin{pf}
	Utilizing the definition of conditional probability, the joint log-likelihood function can be expressed as
	\begin{equation}
		\label{eq3.1.3}
		\setlength{\abovedisplayskip}{3pt}
		\setlength{\belowdisplayskip}{3pt}
		\textstyle
		\begin{aligned}
			\log p_{{\theta}}(X(T),Y(T)) = &\textstyle\log p_{{\theta}}(X(T) | Y(T)) \\&\textstyle+ \log p_{{\theta}}(Y(T)).
		\end{aligned}
	\end{equation}
	Taking conditional expectations on both sides, the following expression is derived
	\begin{equation}
		\label{eq3.1.4}
		\setlength{\abovedisplayskip}{3pt}
		\setlength{\belowdisplayskip}{3pt}
		\textstyle
		\begin{aligned}
			Q({\theta}, \theta_{k}) =&\textstyle L_{{\theta}}(Y(T)) + \int \log p_{{\theta}}(X(T) | Y(T))\cdot\\
			&\textstyle p_{{\theta}_{k}}(X(T) | Y(T)) dX(T). 
		\end{aligned}
	\end{equation}
	Therefore, 
	\begin{equation}
		\label{eq3.1.5}
		\setlength{\abovedisplayskip}{3pt}
		\setlength{\belowdisplayskip}{3pt}
		\textstyle
		\begin{aligned}
			L&\textstyle_{{\theta}}(Y(T))-L_{{\theta}_{k}}(Y(T))= Q({\theta}, \theta_{k}) - Q(\theta_{k}, \theta_{k})\\
			&\textstyle+\int \log\frac{p_{{\theta}_k}(X(T) | Y(T))}{p_{{\theta}}(X(T) | Y(T))}p_{{\theta}_k}(X(T) | Y(T))dX(T)\\
			=&\textstyle Q({\theta}, \theta_{k}) - Q(\theta_{k}, \theta_{k})\\
			&\textstyle+D_{KL}(p_{{\theta}_k}(X(T) | Y(T))||p_{{\theta}}(X(T) | Y(T)))\\
			\geqslant&\textstyle Q(\theta,\theta_k)-Q(\theta_k,\theta_k).
		\end{aligned}
	\end{equation}
	Since the Kullback-Leibler divergence metric satisfies
	\begin{equation}
		\label{eq3.1.6}
		\setlength{\abovedisplayskip}{3pt}
		\setlength{\belowdisplayskip}{3pt}
		\textstyle
		D_{KL}(p_{{\theta}_k}(X(T) | Y(T))||p_{{\theta}}(X(T) | Y(T))) \geqslant 0,\notag
	\end{equation}
	it can be demonstrated that 
	\begin{equation}
		\label{eq3.1.7}
		\setlength{\abovedisplayskip}{3pt}
		\setlength{\belowdisplayskip}{3pt}
		\textstyle
		\begin{aligned}
		L_{\theta_{k+1}}(Y(T))-L_{{\theta}_{k}}(Y(T))& \geqslant \\
		L_{\theta}(Y(T))-L_{\theta_k}(Y(T))&\geqslant\\
		Q({\theta}, \theta_{k}) - Q(\theta_{k}, \theta_{k})&\geqslant 0.
		\end{aligned}
	\end{equation}
	Consequently, maximizing $Q({\theta}, \theta_{k})$ with respect to $\theta$ can result in an increase of the log-likelihood function, i.e., $L_{\theta_{k+1}}(Y(T)) > L_{{\theta}_{k}}(Y(T))$.	{\hfill $\square$}
\end{pf}

\par Utilizing the Bayesian rule and the Markov property of models (\ref{eq2:combined}) and (\ref{eq3.1.1}), it can be inferred that $L_{\theta}(X(T), Y(T))$ satisfies
\begin{equation}
	\label{eq3.1.8}
	\setlength{\abovedisplayskip}{3pt}
	\setlength{\belowdisplayskip}{3pt}
	\textstyle
	\begin{aligned}
		&\textstyle L_{\theta}(X(T), Y(T)) = \log p_{\theta}(Y(T)|X(T)) \\
		&\textstyle\quad+ \log p_{\theta}(X(T))\\
		&\textstyle = \log p_{\theta}(x(1)) \!+\! \sum_{t=1}^{T-1}\log p_{\theta}(x(t+1)|x(t))\\
		&\textstyle\quad+\sum_{t=1}^{T}\log p_{\theta}(y(t)|x(t)).
	\end{aligned}
\end{equation}
Applying the conditional expectation operator $E_{{\theta}_{k}}[\cdot| Y(T)]$ to both sides of (\ref{eq3.1.8}) yields
\begin{equation}
	\label{eq3.1.9}
	\setlength{\abovedisplayskip}{3pt}
	\setlength{\belowdisplayskip}{3pt}
	\textstyle
	Q({\theta}, \theta_{k}) = Q_1+Q_2+Q_3,\tag{11}  
\end{equation}
where,
\begin{subequations}
	\setlength{\abovedisplayskip}{3pt}
	\setlength{\belowdisplayskip}{3pt}
	\begin{align}
		\setlength{\abovedisplayskip}{3pt}
		\setlength{\belowdisplayskip}{3pt}
		\textstyle
		Q_1 =&\textstyle \int \log p_{\theta}(x(1)) p_{\theta_{k}}(x(1)|Y(T)), \label{eq3.1.9a}\\ 
		Q_2 =&\textstyle \sum_{t=1}^{T-1} \iint \log p_{\theta}(x(t+1)|x(t))\cdot\label{eq3.1.9b}\\ 
		&\textstyle p_{\theta_{k}}(x(t+1),x(t)|Y(T))dx(t)dx(t+1), \notag\\
		Q_3 =&\textstyle \sum_{t=1}^{T} \int \log p_{\theta}(y(t)|x(t))\cdot\label{eq3.1.9c}\\
		&\textstyle p_{\theta_{k}}(x(t)|Y(T))dx(t). \notag	
	\end{align}
\end{subequations}
To compute $Q({\theta}, \theta_{k})$, it is necessary to obtain the posterior pdfs of the agent states, such as $p_{\theta_{k}}(x(t)|Y(T))$ and $p_{\theta_{k}}(x(t+1),x(t)|Y(T))$. Due to the uncertain dynamic interactions $g(\cdot)$ and the nonlinear characteristics of the network dynamics, obtaining analytical solutions for the posterior distribution integrals in $Q({\theta},\theta_{k})$ are infeasible. Therefore, a PC-DPF algorithm is employed to yield approximate integrals. This approach also provides a viable numerical computational framework for subsequent parameter estimation.

\section{Distributed State Estimation Using PC-DPF Algorithm}
\label{Sec4}

\par This section introduces a PC-DPF algorithm, aiming to approximate the posterior pdf $p_{\theta_{k}}(x(t)|Y(T))$ using state particles. The calculation relies on the parameter estimate $\theta_k$ obtained from the $k$-th iteration, which will be denoted by the symbol $\theta$ for brevity. Since each agent $v$ can only measure its noisy outputs $y_{v}$ with unknown internal states $x_v$, the proposed algorithm achieves a global consensus state estimation through communication, collaboration, and distributed computation. Fig. \ref{f111} illustrates the local consensus process between agent $v$ and its neighbors, involving multiple interactions with both predecessor neighbors $j$, $j\in\mathcal P_v$ and successor neighbors $j$, $j\in\mathcal S_v$. Individual agents receive particles to obtain global consensus state estimates of the entire network.

\subsection{Sequential Bayesian Estimation}

\par The calculation of the posterior pdf $p_{\theta}(x(t)|Y(T))$ follows the principles of sequential Bayesian estimation. This pdf represents a recursive computation of the confidence level in $x(t)$ given the observations $Y(t)$ up to time $T$, which yields an estimate of the state $x(t)$. It is crucial to address the bilinear coupling between the unknown network dynamic parameter $\theta$ and the latent variable $x(t)$. 

\par According to the models (\ref{eq2:combined}) and (\ref{eq3.1.1}), it can be proven that $p_{{\theta}}(x(t)|Y(t))$ is computable through two sequential steps, using the posterior $p_{{\theta}}(x(t-1)|Y(t-1))$ from the previous time step and the measurements $y(t)$\cite{978374}. In the prediction step, the predicted posterior $p_{{\theta}}(x(t)|Y(t-1))$ is computed as
\begin{equation}
	\label{eq4.1.1}
	\setlength{\abovedisplayskip}{3pt}
	\setlength{\belowdisplayskip}{3pt}
	\textstyle
	\begin{aligned}
		p_{{\theta}}(x(t)|&\textstyle Y(t-1)) = \int p_{{\theta}}(x(t)|x(t-1))\cdot\\ &p_{{\theta}}(x(t-1)|Y(t-1)) dx(t-1).
	\end{aligned}
\end{equation}
Note that (\ref{eq4.1.1}) utilizes the property $p_{\theta}(x(t)|x(t-1),Y(t-1)) = p_{\theta}(x(t)|x(t-1))$ due to the Markov property of the model (\ref{eq3.1.1}). The probabilistic model of state evolution $p_{\theta}(x(t)|x(t-1))$ is defined by (\ref{2a}). In the update step, utilizing the aforementioned Markov property and the recursive Bayesian rule, the predicted posterior $p_{\theta}(x(t)|Y(t-1))$ is updated as
\begin{equation}
	\label{eq4.1.2}
	\setlength{\abovedisplayskip}{3pt}
	\setlength{\belowdisplayskip}{3pt}
	\textstyle
	p_{\theta}(x(t)|Y(t)) =\frac{p_{\theta}({y}(t)|x(t)) p_{\theta}(x(t)|Y(t-1))}{p_{\theta}({y}(t)|Y(t-1))},
\end{equation}
where the normalizing constant
\begin{equation}
	\label{eq4.1.3}
	\setlength{\abovedisplayskip}{3pt}
	\setlength{\belowdisplayskip}{3pt}
	\begin{aligned}
		\textstyle
		p_{\theta}({y}(t)|Y(t-1))=\int& p_{\theta}({y}(t)|x(t))\cdot\\&\textstyle p_{\theta}(x(t)|Y(t-1))dx(t),
	\end{aligned}
\end{equation}
depends on the likelihood function $p_{\theta}({y}(t)|x(t))$ defined by the probabilistic measurement model (\ref{2b}). This process also takes into account the fact that $p_{\theta}( {y}(t)|x(t),Y(t-1)) = p_{\theta}({y}(t)|x(t))$. The recursive computation of $p_{\theta}(x(t)|Y(t))$ starts with $p_{\theta}(x(t)|Y(t))|_{t=0} = p_{\theta}(x_0)$, where the a prior pdf $p_{\theta}(x_0)$ is assumed to be available. Therefore, the estimate of state $x(t)$ can be calculated from the mean of posterior pdf $p_ {\theta} (x (t) | Y (t)) $, which is defined as
\begin{equation}
	\label{eq4.1.4}
	\setlength{\abovedisplayskip}{3pt}
	\setlength{\belowdisplayskip}{3pt}
	\textstyle
	\begin{aligned}
		\hat{x}(t) :=&\textstyle E[x(t)|Y(t)] \\
		=&\textstyle \int x(t) p_{\theta}(x(t)|Y(t)) dx(t).
	\end{aligned}
\end{equation}

\par The recursive propagation of the posterior pdf provides an efficient solution for sequential Bayesian estimation. However, in practical computations, it is often infeasible to compute (\ref{eq4.1.1})-(\ref{eq4.1.4}) involving multi-dimensional integrals. The intricate topology of the large-scale network, represented by $\mathcal{A}$, further complicates the computation of the posterior pdf. Given these challenges, the next subsections introduce the PC-DPF algorithm to yield approximate the distributed computation of the posterior pdf.

\subsection{Particle Smoothing for Individual Agents}
\label{S4.2}
\par To compute the aforementioned posterior pdf, communication within the intricate network topology $\mathcal{A}$ is required. Additionally, the presence of uncertain dynamic interactions $g(x(t)-\mathcal{A}x(t))$ among network agents complicates the accurate calculation of $p_{{\theta}}(x(t)|\mathcal{A}x(t-1))$. Therefore, to reduce the communication requirements and enhance the algorithm's scalability, this subsection first employs the particle smoothing technique to compute the local posterior pdf of individual agents and obtain local state estimation particles. Note that the unknown interaction term $g(x(t)-\mathcal Ax(t))$ will be neglected during the particle smoothing for individual agents, which is inspired by the fact that the term $g(x(t)-\mathcal Ax(t))$ will diminish when reaching the state consensus. To compensate for this operation, the contraction stability will be imposed to yield a robust system identification result.

\par Assume that the state-space model (\ref{1a}) for each agent $v$ at time $t$ consists of a series of particles $\{\tilde{x}_v^{i},i=1,\ldots,M\}$. For a particle $\tilde{x}_v^{i}(t)$, given the local measurements for agent $v$ at time $t$, the conditional probability of the particle is treated as its weight ${w}_v ^{i}(t)$, i.e., $p\{x_v(t) = \tilde{x}_v^{i}(t)|Y_v(t)\} = w_v^{i}(t)$. Then, the posterior pdf of agent $v$ at time $t$ can be approximated as $p_{\theta}(x_v(t)|Y_{v}(t))\approx \sum_{i=1}^{M} w_v^{i}(t)\delta({x}_v(t)-\tilde{x}_v^{i}(t))$. Here, $\delta(\cdot)$ denotes the multi-dimensional Dirac delta function\cite{ZHAO2022137}.

\par According to (\ref{eq4.1.1})-(\ref{eq4.1.4}), the local posterior pdf $p_{\theta}(x_v(t)|Y_{v}(t))$ for agent $v$ can be derived via the prediction step and the update step. In the prediction step, new particles $\tilde{x}_v^{i}(t)$, $i\in\{1,\ldots,M\}$ are sampled from the a priori pdf $p_{\theta}({x}_v(t)|\tilde{x}^{i}_v(t-1))$, according to the probabilistic model (\ref{2a}). These particles are then used to approximate the prediction posterior $p_{\theta}(x_v(t)|Y_{v}(t-1))$, according to
\begin{equation}
	\label{eq4.2.1}
	\setlength{\abovedisplayskip}{3pt}
	\setlength{\belowdisplayskip}{3pt}
	\begin{aligned}
		\textstyle
		p_{\theta}(x_v(t)|&\textstyle Y_{v}(t-1))\approx \sum_{i=1}^{M}{w}_v^{i}(t-1)\cdot\\&\textstyle p_{\theta}({x}_v(t)|\tilde{x}^{i}_v(t-1))\delta({x}_v(t)-\tilde{x}_v^{i}(t)).
	\end{aligned}
\end{equation}
In the update step, the weights $\{{w}_v^{i}(t)\}_{i=1}^M$ are updated based on the new local observations $y_{v}(t)$ from agent $v$ with the predicted particles obtained in the prediction step as
\begin{equation}
	\label{eq4.2.2}
	\setlength{\abovedisplayskip}{3pt}
	\setlength{\belowdisplayskip}{3pt}
	\textstyle
	w_v^{i}(t) = \frac{p_{\theta}(y_{v}(t) | \tilde{x}_v^{i}(t))}{\sum_{j=1}^{M}p_{\theta}(y_{v}(t) | \tilde{x}_v^{j}(t))}.
\end{equation}
Subsequently, the particles $\{\tilde{x}_v^{i}(t)\}_{i=1}^M$ are resampled with weights $\{{w}_v^{i}(t)\}_{i=1}^M$ to obtain updated particles $\{\tilde{x}_v^{j}(t+1)\}_{j=1}^M$ at time $t+1$. Repeat (\ref{eq4.2.1})-(\ref{eq4.2.2}) until $t=T$ and return the state particles and weights $\{\tilde{x}_v^{i}(t),w_v^{i}(t)\}_{i=1}^M$, for $t=\{1,\ldots,T\}$. The recursive computation of $p_{\theta}(x_v(t)|Y_v(t))$ starts with the initialized particles $\{\tilde{x}_{v}^{i}(1)\}_{i=1}^M\sim  p_{\theta}({x}_{v}(1))$.

\par To address the problem of inaccurate state approximation in continuous systems caused by long discrete time intervals for discretization, a smoothing estimate is provided for the entire state trajectory $x_v(1:T)$. This involves extracting the necessary state subsequences and simultaneously adjusting the weights during the backward pass, yielding weighted samples aimed at smoothing the distribution at fixed time intervals. Based on the law of total probability and the Bayesian rule, it follows that \cite{10.1214/14-STS511}
\begin{equation}
	\label{eq4.2.3}
	\setlength{\abovedisplayskip}{3pt}
	\setlength{\belowdisplayskip}{3pt}
	\textstyle
	\begin{aligned}
		&\quad p_{\theta}(x_v(t)|Y_{v}(T))= p_{\theta}(x_v(t)|Y_{v}(t))\cdot\\\textstyle\int&\textstyle\frac{p_{\theta}(x_v(t+1)|x_v(t))p_{\theta}(x_v(t+1)|Y_{v}(T))}{p_{\theta}(x_v(t+1)|Y_{v}(t))}dx_v(t+1).
	\end{aligned}
\end{equation} 
Given the local observations $Y_{v}(T)$, the smoothing distribution $p_{\theta}(x_v(t) | Y_{v}(T))$ of the states can be approximated as
\begin{equation}
	\label{eq4.2.4}
	\setlength{\abovedisplayskip}{3pt}
	\setlength{\belowdisplayskip}{3pt}
	\textstyle
	p_{\theta}(x_v(t)|Y_{v}(T))\approx \sum_{i=1}^{M} w_v^{i}(t|T)\delta(x_v(t)-\tilde{x}_v^{i}(t)),
\end{equation} 
where the backward update weight $w_v^{i}(t|T)$ of the $i$-th particle can be deduced from (\ref{eq4.2.1})-(\ref{eq4.2.3}) as
\begin{equation}
	\label{eq4.2.5}
	\setlength{\abovedisplayskip}{3pt}
	\setlength{\belowdisplayskip}{3pt}
	\begin{aligned}
		\textstyle
		w_v^{i}(t|T)=\sum_{j=1}^{M}&\textstyle w_v^{j}(t+1|T)\cdot\\&\textstyle\frac{w_v^{i}(t)p_{\theta}(\tilde{x}_v^{j}(t+1)|\tilde{x}_v^{i}(t))}{\sum_{k=1}^{M}w_v^{k}(t)p_{\theta}(\tilde{x}_v^{j}(t+1)|\tilde{x}_v^{k}(t))}.
	\end{aligned}
\end{equation}
The derivation of this weight involves the application of Bayes' rule to the posterior $p_{\theta}(x_v(t+1),x_v(t)|Y_v(T))$, which satisfies
\begin{equation}
	\label{eq4.2.6}
	\setlength{\abovedisplayskip}{3pt}
	\setlength{\belowdisplayskip}{3pt}
	\begin{aligned}
		\textstyle
		&p_{\theta}(x_v(t+1),x_v(t)|Y_v(T))\\
		=&\textstyle
		p_{\theta}(x_v(t)|x_v(t+1),Y_v(t))
		p_{\theta}(x_v(t\!+\!1)|Y_v(T))\\
		=&\textstyle\frac{p_{\theta}(x_v(t+1)|x_v(t))p_{\theta}(x_v(t)|Y_v(t))}{p_{\theta}(x_v(t+1)|Y_v(t))}p_{\theta}(x_v(t\!+\!1)|Y_v(T)).
	\end{aligned}
\end{equation}
A detailed derivation can be found in Lemma 6.1 of \cite{SCHON201139}.

\par Resample the particles $\{\tilde{x}_v^{i}(t)\}_{i=1}^M$ according to the smoothed weights $\{w_v^{i}(t|T)\}_{i=1}^\frac{M}{{J_v}_{\text{max}}}$ to find the $\frac{M}{{J_v}_{\text{max}}}$ particles $\{\tilde{x}_v^{j}(t)\}_{j=1}^{\frac{M}{{J_v}_{\text{max}}}}$ with the largest weights, and compute the set of state estimation particles corresponding to each agent $v$, satisfying $\hat{x}_v^{j}(t)=\tilde{x}_v^{j}(t)\cdot w_v^{j}(t|T)$, $j=\{1,\ldots,\frac{M}{{J_v}_{\text{max}}}\}$. Here, ${J_v}_{\text{max}}$ is chosen as the maximum number of neighbors among all agents in the network, aiming to improve the efficiency of data transmission. These state estimation particles $\{\hat{x}_v^{j}(t)\}_{j=1}^{\frac{M}{{J_v}_{\text{max}}}}$ for each agent $v$ will be combined with the PC-DPF algorithm proposed in the next subsection to realize global consensus state estimation of latent states in large-scale networks. 

\subsection{The Proposed PC-DPF Algorithm}
\label{4.3}
\par This subsection introduces a global consensus algorithm based on local state estimation particles of each to facilitate the distributed computation and information exchange among neighboring agents.

\par In each communication round, agents randomly transmit a subset of representative particles to their neighbors, following the known directed network topology $\mathcal{A}$. This dynamic information propagation, analogous to the spread of an epidemic, ensures high fault tolerance and self-stabilization \cite{10.1145/41840.41841}. The communication incurs moderate overhead and benefits from simplicity, as it operates without the need for error recovery mechanisms.

\par In Section \ref{S4.2}, each agent generates local smoothing state estimation particles using its own noisy observations. However, for the considered isomorphic network dynamics (\ref {eq3.1.1}), which involves unknown interactions ${g}(x(t) - \mathcal{A}x(t))$ among agents, and where each agent has unique local process and measurement noises, it becomes necessary to diffuse local state estimation particles across the entire networked system.

\par To address the unknown interactions in the network, this study tackles the following particle aggregation problem: in a network composed of $V$ agents, each holding a set of particles $\{\hat{x}_v^{i}(t)\}_{i=1}^{\frac{M}{{J_v}_{\text{max}}}}$, the goal is to compute an aggregation function (e.g., mean) for all agents' particles in a fully decentralized and fault-tolerant manner. The proposed computation method converges exponentially to consensus particles for the entire network, as demonstrated by Theorem \ref{thm3}. As a result, the unknown interaction term $g(x(t)-\mathcal Ax(t))$ will diminish. Therefore, the omission of this term in the smoothing step will not affect the final convergence result.

\begin{algorithm}
	\caption{Particle Consensus-based DPF algorithm (PC-DPF)}
	\label{al3}
	\begin{algorithmic}
		\State \hspace*{-1em}1. Initialize the iteration number $I_{con}$, and the number of neighbors $J_v$ for each agent $v$.
		\State \hspace*{-1em}2. At each iteration time $t$, perform the following steps:
		\State \hspace*{1em}2.1 Let $Xc_{t,v}=\sum_j Xc_{t-1,j}$, $n_{t,v}=\sum_j n_{t-1,j}$, where $j\in\mathcal P_v$. 
		\State \hspace*{1em}2.2 Broadcast $(\frac{1}{2}Xc_{t,v},\frac{1}{2}n_{t,v})$ to a randomly selected neighbor $j\in\mathcal P_v$ from its total neighbors $J_v$ and itself $v$.
		\State \hspace*{1em}2.3 Compute the estimate of the average $\frac{Xc_{t,v}}{n_{t,v}}$ at step $t$ and the convergence error $e$:
		\begin{equation} 
			\label{eq16}
			\setlength{\abovedisplayskip}{3pt}
			\setlength{\belowdisplayskip}{3pt}
			\textstyle
			e=\text{mean}(\{Xc_{t,v}\}_{v=1}^{V})-\frac{Xc_{t,v}}{n_{t,v}}.
		\end{equation}
		\State \hspace*{-1em}3. Terminate the iteration when the iteration number reaches $I_{con}$ or the convergence error $e\leqslant1.5\times\delta$\footnotemark{}.
		\State \hspace*{-1em}4. Generate the global state estimation particle set  $\{\bar{\hat{x}}^{i}(t)\}_{i=1}^{\lceil \frac{M}{{J_v}_{\text{max}}}\rceil \times V}$ of each agent $v$ after communication utilizing (\ref{eq4.3.1}).
		\State \hspace*{-1em}5. Repeat steps 2-4 until $t=T$.
		\State \hspace*{-1em}6. Return the global consensus state estimation particles 	$\{\bar{\hat{x}}^{i}(t)\}_{i=1}^{\lceil \frac{M}{{J_v}_{\text{max}}}\rceil \times V}$, for $t=\{1,\ldots,T\}$.
	\end{algorithmic}
\end{algorithm}
\footnotetext{$\delta$ is the error derived in Theorem \ref{thm3}, which can be chosen as a user-defined constant during the actual algorithm execution. Since this error value is small, $1.5\times\delta$ is used here to represent the convergence error condition for iteration termination only to show that this iteration termination condition can be larger than $\delta$ and to avoid abuse of notation.}

\par Specifically, each agent $v$ initializes a particle transmission vector $Xc_{0,v}\in\mathbb{R}_{\lceil\frac{M}{{J_v}_{\text{max}}}\rceil\times V}$, as $Xc_{0,v}((v-1)\lceil\frac{M}{{J_v}_{\text{max}}}\rceil+1:v\lceil\frac{M}{{J_v}_{\text{max}}}\rceil)=\{\hat{x}_v^{i}(0)\}_{i=1}^{\frac{M}{{J_v}_{\text{max}}}}$. Simultaneously, each agent $v$ maintains a propagation counter with an initial value of $n_{0,v}=1$. As outlined in Algorithm \ref{al3}, at each time step $t$, each agent $v$ aggregates the particle transmission vectors and counters from its predecessor neighbors $j$, $j\in\mathcal P_v$ and updates its local particle transmission vector and counter as $(\frac{1}{2}Xc_{t,v},\frac{1}{2}n_{t,v})$. The agent $v$ then randomly selects a successor neighbor $j$, $j\in\mathcal S_v$, and transmits $(\frac{1}{2}Xc_{t,v},\frac{1}{2}n_{t,v})$ to it. Theorem \ref{thm3} establishes that when $t=\log_2 V+\log_2 \frac{1}{\bar \delta}$, all agents in the network, with a probability not less than $1-\frac{1}{V}$, possess a consensus particle transmission vector $Xc_{t,v}$ with a relative error not exceeding $\delta$ compared to the true value $\text{mean}(\{Xc_{0,v}\}_{v=1}^{V})$. Then the global consensus state estimation particles of all agents in the network can be expressed as 		
\begin{equation}
	\label{eq4.3.1}
	\setlength{\abovedisplayskip}{3pt}
	\setlength{\belowdisplayskip}{3pt}
	\textstyle
	\{\bar{\hat{x}}^{i}(t)\}_{i=1}^{\lceil \frac{M}{{J_v}_{\text{max}}}\rceil \times V}=V \times \frac{Xc_{t,v}}{n_{t,v}}.
\end{equation}
Here, multiplying the right-hand side by $V$ is necessary, because after achieving consensus among $V$ agents on $\frac{Xc_{t,v}}{n_{t,v}}$, the magnitude of this value is reduced to $\frac{1}{V}$. Thus, the multiplication by $V$ is applied to restore it to its original magnitude. The detailed flow of the algorithm is shown in Algorithm \ref{al3}. The algorithm ultimately generates the state estimation particles $\{\bar{\hat{x}}^{i}(t)\}_{i=1}^{\lceil \frac{M}{{J_v}_{\text{max}}}\rceil \times V}$ that contain all the state features of all agents in the network.

\par In summary, the proposed particle aggregation method can achieve an approximate consensus for all agents in the network at time $T$. Due to the random selection of successor neighbors for particle broadcasting, the algorithm avoids bottlenecks or single points of failure. Furthermore, the algorithm exhibits robustness in the face of dynamic network topology and unreliable network conditions, making it well-suited for applications involving uncertain dynamic interactions. The application of this algorithm facilitates the calculation of posterior pdfs $p_{\theta_{k}}(x(t)|Y(T))$ and $p_{\theta_{k}}(x(t+1),x(t)|Y(T))$ in networks with dynamic interactions, which will be detailed in Section \ref{Sec5}.

\section{Distributed Dynamic Parameter Identification Based on EM Framework}
\label{Sec5}

\par In this section, a distributed parameter identification method under the EM framework will be presented. The method is designed for globally incrementally $l_2$ contraction stable nonlinear state space models in networks. Based on the distributed state estimation discussed in Section \ref{Sec4}, the proposed identification method will address the bilinear coupling between the unknown parameters and unmeasurable states. A crucial aspect of this approach lies in its ability to ensure the contraction stability of the large-scale network, a key factor for reliable network identification.

\subsection{E-step: Distributed Calculation of the  Likelihood Function Approximation $Q({\theta}, \theta_{k})$}

\par The sequential Bayesian estimation method discussed in Section \ref{Sec4} is now employed to calculate the approximation of $Q_1$, $Q_2$, and $Q_3$ in (\ref{eq3.1.9}). Based on the definitions in (\ref{eq4.2.3})-(\ref{eq4.2.5}), the particle form of $Q ({\theta}, \theta_{k}) $ can be approximated as
\begin{equation}
	\label{eq5.1.1}
	\setlength{\abovedisplayskip}{3pt}
	\setlength{\belowdisplayskip}{3pt}
	\textstyle
	\tilde{Q}({\theta}, \theta_{k}) = 	\tilde{Q}_1+\tilde{Q}_2+\tilde{Q}_3,\tag{24}  
\end{equation}
where,
\begin{subequations}
	\setlength{\abovedisplayskip}{3pt}
	\setlength{\belowdisplayskip}{3pt}
	\begin{align}
		\setlength{\abovedisplayskip}{3pt}
		\setlength{\belowdisplayskip}{3pt}
		\textstyle
		\tilde{Q}_1 :=&\textstyle \sum_{i=1}^{M} w_v^{i}(1|T) 
		\log p_{\theta}(\tilde{x}_v^i(1)),\label{eq5.1.1a}\\ 
		\tilde{Q}_2 :=&\textstyle\sum_{t=1}^{T-1}\sum_{i=1}^{M}
		\sum_{j=1}^{M} w_v^{ij}(t|T)\cdot
		\label{eq5.1.1b}\\
		&\qquad\qquad\log p_{\theta}(\tilde{x}_v^j(t+1)|\tilde{x}_{\mathcal P_v}^i(t)),
		\notag\\ 
		\tilde{Q}_3 :=&\textstyle \sum_{t=1}^{T}
		\sum_{i=1}^{M}w_v^{i}(t|T)\cdot		\label{eq5.1.1c} \\
		&\qquad\qquad\log p_{\theta}(y_v(t)|\tilde{x}_{\mathcal P_v}^i(t)).\notag
	\end{align}
\end{subequations}
The weight $w_v^{ij}(t|T)$ is given by
\begin{equation}
	\label{eq5.1.2}
	\setlength{\abovedisplayskip}{3pt}
	\setlength{\belowdisplayskip}{3pt}
	\textstyle
	w_v^{ij}(t|T) = \frac{w_t^i w_{t+1|T}^jp_{\theta_{k}}(\tilde{x}_v^j(t+1)|\tilde{x}_{\mathcal P_v}^i(t))}{\sum_{l=1}^{M}w_t^l p_{\theta_{k}}(\tilde{x}_v^j(t+1)|\tilde{x}_{\mathcal P_v}^l(t))}.
\end{equation}

\par Note that the computation of $\tilde{Q}({\theta}, \theta_{k})$ needs to consider the dynamic influence of each agent $v$'s predecessor neighbors ${\mathcal P_v}$. Therefore it is necessary to calculate $p_{\theta}(\tilde{x}_v^j(t+1)|\tilde{x}_{\mathcal P_v}^i(t))$ and $p_{\theta}(y_v(t)|\tilde{x}_{\mathcal P_v}^i(t))$ based on the model (\ref{eq1:combined}). Due to the uncertainty of ${g}(x(t)-\mathcal{A}x(t))$, these two pdfs remain difficult to compute accurately. To address this challenge, Algorithm \ref{al3} assigns global consensus state estimation particles to individual agents, achieving efficient distributed computation of $\tilde{Q}({\theta}, \theta_{k})$.

\par The global consensus state estimation particle $\bar{\hat{x}}^i(t)$ cannot independently obtain the corresponding weight $\bar{w}^i(t)$ and particle value $\bar{\tilde{x}}^i(t)$ similar to (\ref{eq5.1.1}). By assuming the i.i.d noises $\varepsilon_v$ and $\eta_v$ to follow Gaussian distribution, the function $\log p_{\theta}(\cdot)$ is concave, as $\log$ is concave and monotonically increasing, and $p_{\theta}$ is concave \cite{kkt}. According to the Jensen's inequality, it can be inferred that 
\begin{equation}
	\setlength{\abovedisplayskip}{3pt}
	\setlength{\belowdisplayskip}{3pt}
	\begin{aligned}
		\textstyle
		\sum_{i=1}^{\lceil \frac{M}{{J_v}_{\text{max}}}\rceil \times V} &\bar{w}^i(t) \log p_{\theta}(\cdot| \bar{\tilde{x}}^i(t))\leqslant\\
		&\textstyle
		\log p_{\theta_v}(\cdot| \sum_{i=1}^{\lceil \frac{M}{{J_v}_{\text{max}}}\rceil \times V}\bar{w}^i(t|T)\bar{\tilde{x}}^i(t)).\notag
	\end{aligned}
\end{equation}
Based on the definition in (\ref{eq4.1.4}), it is known that $\sum_{i=1}^{\lceil \frac{M}{{J_v}_{\text{max}}}\rceil \times V}\bar{w}^i(t|T)\bar{\tilde{x}}^i(t) =\sum_{i=1}^{\lceil \frac{M}{{J_v}_{\text{max}}}\rceil \times V}\bar{\hat{x}}^i(t)=\bar{\hat{x}}(t)$. Therefore, $\sum_{i=1}^{\lceil \frac{M}{{J_v}_{\text{max}}}\rceil \times V} \bar{w}^i(t|T)\log p_{\theta}(\cdot| \bar{\tilde{x}}^i(t))$ can be directly substituted with $\log p_{\theta}(\cdot|\bar{\hat{x}}(t))$.

\par By employing particles that exhibit global consistency, the global consensus particle approximate likelihood function $\bar{\hat{Q}}({\theta}, \theta_{k})$ can be defined as
\begin{equation}
	\label{eq5.1.4}
	\setlength{\abovedisplayskip}{3pt}
	\setlength{\belowdisplayskip}{3pt}
	\textstyle
	\bar{\hat{Q}}({\theta}, \theta_{k})=\bar{\hat{Q}}_1+\bar{\hat{Q}}_2+\bar{\hat{Q}}_3,\tag{26}  
\end{equation}
where,
\begin{subequations}\label{combined:eq4.1.5}
	\setlength{\abovedisplayskip}{3pt}
	\setlength{\belowdisplayskip}{3pt}
	\begin{align}
		\setlength{\abovedisplayskip}{3pt}
		\setlength{\belowdisplayskip}{3pt}
		\textstyle
		\bar{\hat{Q}}_1 :=&\textstyle 
		\log p_{\theta}(\bar{\hat{x}}(1)),\label{eq4.1.5a}\\ 
		\bar{\hat{Q}}_2 :=&\textstyle \sum_{t=1}^{T-1}\log p_{\theta}(\bar{\hat{x}}(t+1)|\bar{\hat{x}}(t)),
		\label{eq4.1.5b}\\ 
		\bar{\hat{Q}}_3 :=&\textstyle \sum_{t=1}^{T} 
		\log p_{\theta}({y}_v(t)|\bar{\hat{x}}(t)).
		\label{eq4.1.5c} 	
	\end{align}
\end{subequations}
This formula represents the distributed computation function for each agent using the global consensus state estimation particles. As each agent can perform distributed computation, the identification of unknown parameters can be carried out at a single agent in an isomorphic network.

\subsection{M-Step: Parameter Identification with  Stability Constraint} 

\par To ensure accurate identification of model parameters in the presence of uncertain dynamic interactions ${g}(x(t) - \mathcal{A}x(t))$, this study extends a set of stabilizing models constructed based on the deterministic portion of the model (\ref{1a}), defined as
\begin{equation}
	\label{5.2.1}
	\setlength{\abovedisplayskip}{3pt}
	\setlength{\belowdisplayskip}{3pt}
	\textstyle
	{x}_v(t+1) =\sum_{i=0}^{q}{{\theta}_{i}}{{f}_{i}}(x_v(t), \varepsilon_v(t)).
\end{equation}
By linearizing around the initial conditions, the differential dynamics of (\ref{5.2.1}) can be expressed as
\begin{equation}
	\label{5.2.2}
	\setlength{\abovedisplayskip}{3pt}
	\setlength{\belowdisplayskip}{3pt}
	\textstyle		
	{{\Delta }_v}(t+1)=F(x_v (t),{\theta})\Delta_v (t), 
\end{equation} 
where, $F(x_v,{\theta})=\sum_{i=0}^{q}{{\theta}_{i}}\frac{\partial}{\partial x_v}{{f}_{i}}(x_v, \varepsilon_v)$. For all nontrivial parameterizations of the above model, explicit and convex properties of the set $\Theta$ are defined to represent all stable model parameters. The contraction stability of an agent model is guaranteed by the following theorem.

\begin{thm}
	\label{thm4}
	If there exists a symmetric matrix $P_v>0$ and $\kappa>0$, such that
	\begin{equation}
		\label{5.2.3}
		\setlength{\abovedisplayskip}{3pt}
		\setlength{\belowdisplayskip}{3pt}
		\textstyle
			\Theta:=\left\lbrace \theta\in\begin{bmatrix}
				(2-\kappa)I_m-P_v & {F}(x_v,\theta)^\mathrm{T}\\
				{F}(x_v,\theta) & P_v
			\end{bmatrix}\succeq 0\right\rbrace.
	\end{equation}
	Then the model (\ref{5.2.1}) is said to be globally incrementally $l_2$ contraction stable.
\end{thm}
\begin{pf}
	Applying the Schur complementary lemma \cite{sc} to (\ref{5.2.3}), translates into
	\begin{equation}
		\label{5.2.4}
		\setlength{\abovedisplayskip}{3pt}
		\setlength{\belowdisplayskip}{3pt}
		\textstyle
		(2-\kappa)I_m-P_v-{F}(x_v,\theta)^\mathrm{T}{P}^{-1}_v{F}(x_v,\theta) I\geqslant0
	\end{equation}
	Since $|I_m-P_v|^2_{{P}^{-1}_v}={P}^{-1}_v-2I_m+P_v\geqslant0$, it follows that $ 2I_m-P_v \leqslant {P}^{-1}_v$. Substituting into (\ref{5.2.4}), it has
	\begin{equation}
		\label{5.2.5}
		\setlength{\abovedisplayskip}{3pt}
		\setlength{\belowdisplayskip}{3pt}
		\textstyle
		{F}(x_v,\theta)^\mathrm{T}{P}^{-1}_v{F}(x_v,\theta)- {P_v}^{-1}\leqslant-\kappa I_m
	\end{equation}
Simultaneously multiplying ${{\Delta }_v}^\mathrm{T}$ on the left and ${{\Delta }_v}$ on the right for both sides of (\ref{5.2.5}), it can be derived that
\begin{equation}
	\label{5.2.6}
	\setlength{\abovedisplayskip}{3pt}
	\setlength{\belowdisplayskip}{3pt}
	\textstyle
		|F(x_v,{\theta})\Delta_v |^2_{{P}^{-1}_v}-|\Delta_v| ^2_{{P}^{-1}_v}\leqslant-\kappa|\Delta_v| ^2_{{P}^{-1}_v}
\end{equation}
That is, the differential storage function $V=|\Delta_v| ^2_{{P}^{-1}_v}$ can be found to ensure that the model (\ref{5.2.1}) is globally incrementally $l_2$ contraction stabilized under the stabilization constraint defined in (\ref{5.2.3}), as described in Definition \ref{def1}.
{\hfill $\square$}
\end{pf}
\begin{rem}
	\label{rem}
	It can be seen that the stability constraint can be expressed in quadratic form as shown in (\ref{5.2.5}) and $x_v$ can be replaced by the global particle state estimations obtained in Section \ref{4.3}. As a result, the contraction stability constraint is convex with respect to $\theta$.
\end{rem}

\par With an approximation $\bar{\hat{Q}}({\theta}, \theta_{k})$ of the function ${Q}({\theta}, \theta_{k})$, attention now turns to the M-step calculation. To ensure that the nonlinear state-space model obtained from the identification remains stable, this study combines the M-step with the stability constraint defined in (\ref{5.2.3}). In this step, the goal is to maximize the constrained $\bar{\hat{Q}}({\theta}, \theta_{k})$ with respect to ${\theta}$. Consequently, the maximum likelihood estimate can be expressed as the following optimization problem as
\begin{equation}
	\label{eq5.2.7}
	\setlength{\abovedisplayskip}{3pt}
	\setlength{\belowdisplayskip}{3pt}
	\textstyle
	\hat{\theta}=\underset{\theta\in\Theta}{\arg \max} \;\bar{\hat{Q}}({\theta}, \theta_{k}).  
\end{equation}
For i.i.d noises $\varepsilon_v$ and $\eta_v$ to follow Gaussian distribution, the expression (\ref{eq5.1.4}) can be equivalently expressed as
\begin{equation}
	\label{eq5.2.8}
	\setlength{\abovedisplayskip}{3pt}
	\setlength{\belowdisplayskip}{3pt}
	\begin{aligned}
	\textstyle
	\bar{\hat{Q}}({\theta}, \theta_{k})=&\textstyle-\sum_{t=1}^{T-1}|\bar{\hat{x}}(t+1)-\theta f(\bar{\hat{x}}(t))-g(0)|^2\\
	&-\textstyle\sum_{t=1}^{T}
	|y_v(t)- \theta h(\bar{\hat{x}}(t))|^2\\&-(2T-1)\theta\log |\theta|, v\in\mathcal{V}.
	\end{aligned}
\end{equation}
As stated in Remark \ref{rem}, once the global particle state estimate $\bar{\hat{x}}(t)$ in $\bar{\hat{Q}}({\theta}, \theta_{k})$ has been determined, it is easy to show that $\frac{\partial^2}{\partial \theta^2}\bar{\hat{Q}}({\theta}, \theta_{k})\leqslant0$, i.e. the objective function is convex. Note that the mean and variance of the noise $\varepsilon_v$ and $\eta_v$ have been incorporated into the parameter $\theta$ for identification, and hence the existence of $|\theta|$ guarantees that the noise parameter must be greater than zero.

\begin{thm}
	\label{thm5}
	For the optimization problem (\ref{eq5.2.7}), there exists at least one optimal solution $\theta^*$ such that
	\begin{equation}
	\theta^* =\underset{\theta\in\Theta}{\arg \max} \;\bar{\hat{Q}}({\theta}, \hat{\theta}_{k}). 
	\end{equation}
\end{thm}
\begin{pf}
	By introducing a non-positive Lagrange multiplier $\lambda\leqslant0$, the problem (\ref{eq5.2.7}) can be transformed into an unconstrained optimization problem, yielding the Lagrange function
	\begin{equation}
		\label{eq5.2.9}
		\setlength{\abovedisplayskip}{3pt}
		\setlength{\belowdisplayskip}{3pt}
		\begin{aligned}
			\textstyle
			L(\theta,\lambda)=-&\textstyle\sum_{t=1}^{T-1}|\bar{\hat{x}}(t+1)-\theta
			f(\bar{\hat{x}}(t))
			-g(0)|^2\\
			&-\textstyle\sum_{t=1}^{T}
			|y_v(t)- \theta h(\bar{\hat{x}}(t))|^2\\
			&-(2T-1)\theta\log |\theta|\\
			&\textstyle+\lambda(|\theta
			\frac{\partial}{\partial \bar{\hat{x}}}f(\bar{\hat{x}}(t))|^2_{{P}^{-1}_v}-{{P}^{-1}_v}+\kappa I_m).  
		\end{aligned}
	\end{equation}
	Therefore, according to the necessary conditions for an optimal solution (extreme value condition), it can be obtained that
	\begin{equation}
		\label{eq5.2.10}
		\setlength{\abovedisplayskip}{3pt}
		\setlength{\belowdisplayskip}{3pt}
		\begin{aligned}
			\textstyle
			\frac{\partial}{\partial \theta}&L(\theta,\lambda)=\\&\textstyle\sum_{t=1}^{T-1}2(\bar{\hat{x}}(t+1)-\theta
			f(\bar{\hat{x}}(t))
			-g(0))f(\bar{\hat{x}}(t))\\
			&+\textstyle\sum_{t=1}^{T}2
			(y_v(t)-\theta h(\bar{\hat{x}}(t)))h(\bar{\hat{x}}(t))\\
			&-(2T-1)\log |\theta|-(2T-1)\\
			&\textstyle+2\lambda{{P}^{-1}_v}(\frac{\partial}{\partial \bar{\hat{x}}}f(\bar{\hat{x}}(t)))^2\theta.  
		\end{aligned}
	\end{equation}
	Since
		\begin{equation}
		\label{eq5.2.10.2}
		\setlength{\abovedisplayskip}{3pt}
		\setlength{\belowdisplayskip}{3pt}
		\begin{aligned}
			\textstyle
			&\textstyle\frac{\partial^2}{\partial \theta^2} L(\theta,\lambda)=-\textstyle\frac{2T-1}{|\theta|}-\sum_{t=1}^{T-1}2(
			f(\bar{\hat{x}}(t))
			)^2\\
			&-\textstyle\sum_{t=1}^{T}2
			(h(\bar{\hat{x}}(t)))^2
			\textstyle+2\lambda{{P}^{-1}_v}(\frac{\partial}{\partial \bar{\hat{x}}}f(\bar{\hat{x}}(t)))^2\leqslant0,  
		\end{aligned}
	\end{equation}
	the function $\frac{\partial}{\partial \theta}L(\theta, \lambda)$ is monotonically decreasing. As $\theta\to0$, $\frac{\partial}{\partial \theta}L(\theta,\lambda)\to\infty$, and as $\theta\to\infty$, $\frac{\partial}{\partial \theta}L(\theta,\lambda)\to-\infty$, there must exist $\theta^*$ such that $\frac{\partial}{\partial \theta}L(\theta,\lambda)|_{ \theta=\theta ^*}=0$.
	
	\par Moreover, when $\lambda\neq0$, a set of $\kappa$ and $P_v$ can be found such that $\theta ^*$ satisfies
		\begin{equation}
		\label{eq5.2.11}
		\setlength{\abovedisplayskip}{3pt}
		\setlength{\belowdisplayskip}{3pt}
		\textstyle
		\theta ^*=\frac{\sqrt{ I_m -\kappa I_m{{P}_v} }}{\frac{\partial}{\partial \bar{\hat{x}}}f(\bar{\hat{x}}(t))},
	\end{equation}
	which implies that
	\begin{equation}
		\label{eq5.2.11.2}
		\setlength{\abovedisplayskip}{3pt}
		\setlength{\belowdisplayskip}{3pt}
		\textstyle
		\frac{\partial}{\partial \lambda}L(\theta^*,\lambda)=|\theta^*
		\frac{\partial}{\partial \bar{\hat{x}}}f(\bar{\hat{x}}(t))|^2_{{P}^{-1}_v}-{{P}^{-1}_v}+\kappa I_m=0,  
	\end{equation}
	and
		\begin{equation}
		\label{eq5.2.11.3}
		\setlength{\abovedisplayskip}{3pt}
		\setlength{\belowdisplayskip}{3pt}
		\textstyle
		\lambda(|\theta^*
		\frac{\partial}{\partial \bar{\hat{x}}}f(\bar{\hat{x}}(t))|^2_{{P}^{-1}_v}-{{P}^{-1}_v}+\kappa I_m)=0.  
	\end{equation}
	
	\par Thus, this ensures that at least one extreme point $\theta ^*$ satisfies the Karush-Kuhn-Tucker (KKT) conditions \cite{kkt}. Since the optimization problem is convex, this extreme point is the optimum point.
	{\hfill $\square$}
\end{pf}

\section{Particle Consensus-based Distributed Particle EM algorithm}
\label{Sec6}

\par This section proposes a PC-DPEM algorithm for the identification of model parameters in large-scale isomorphic nonlinear networks, which firstly describes the key steps of the algorithm and follows with a comprehensive analysis of its performance.

\subsection{The Proposed PC-DPEM Algorithm}

\begin{algorithm}
	\caption{PC-DPEM algorithm}
	\label{al4}
	\begin{algorithmic}
		\State \hspace*{-1em}${1.}$ Set $k=0$, and initialize $\theta_{k}$ such that $L_{\theta_{k}}(Y)$ is finite.
		\State \hspace*{-1em}2. $\mathbf{E-step\ Calculation}$:
		\State \hspace*{0.8em}2.1 Generate
		\begin{equation}
		\setlength{\abovedisplayskip}{3pt}
		\setlength{\belowdisplayskip}{3pt}
		\{\bar{\hat{x}}^{i}(t)\}_{i=1}^{\lceil \frac{M}{{J_v}_{\text{max}}}\rceil \times V}=\textbf{PC-DPF}(Xc_{0,v},w_{0,v}).\notag
		\end{equation}
		\State \hspace*{0.8em}2.2 Use $\{\bar{\hat{x}}(t)\}$ to calculate $\bar{\hat{Q}}({\theta}, \theta_{k})$ in (\ref{eq5.1.4}).
		\State \hspace*{-1em}3. $\mathbf{M-step\ Calculation}$:
		\State \hspace{0.8em}Solve the optimization problem in (\ref{eq5.2.7}).
		\State \hspace*{-1em}4. Check the following stopping criterion
		\begin{equation}
			\setlength{\abovedisplayskip}{3pt}
			\setlength{\belowdisplayskip}{3pt}
			\textstyle
			\bar{\hat{Q}} ({\theta}_{k+1}, {\theta}_{k})-\bar{\hat{Q}} ({\theta}_{k}, {\theta}_{k}) < \epsilon \notag
		\end{equation}
		for some user-chosen $\epsilon > 0$\footnotemark{}. If not satisfied, update $k = k + 1$ and return to step 2; otherwise, terminate.  
	\end{algorithmic}
\end{algorithm}
\footnotetext{$\epsilon$ is a user-defined constant that serves as an iterative termination condition for Algorithm \ref{al4}. The definition of this value depends on the parameter solution error allowed by the user.}

\par Note that this algorithm is sensitive to the number of particles. A larger $M$ results in more accurate identification results. However, it also increases the computational complexity. Therefore, it is necessary to strike a balance between the identification accuracy and computational complexity.

\subsection{Communication Overhead and Computational Complexity}
\label{6.2}

\par The transmission of particles across the network incurs communication overhead. In Algorithm \ref{al3}, each agent $v$ broadcasts $n\lceil \frac{M}{{J_v}_{\text{max}}}\rceil\times V+1$ scalar values in parallel, where $n$ is the dimension of the state space. For the network topology considered in this study, the number of agent $v$'s neighbors typically exhibits a complexity of $\mathcal{O}({J_v}_{\text{max}})$. Therefore, the worst-case communication complexity for each agent in each iteration is $\mathcal{O}((n\lceil \frac{M}{{J_v}_{\text{max}}}\rceil\times V+1){{J_v}_{\text{max}}})$. Since all agents in the network are allowed to broadcast in parallel to any one of their successor neighbors, the number of agents making broadcasts in the network at a single iteration is determined by the number of directed edges $E$ in the communication topology. Therefore, the worst-case communication complexity of the entire network at each iteration is $\mathcal{O}(E(n\lceil \frac{M}{{J_v}_{\text{max}}}\rceil\times V+1){{J_v}_{\text{max}}})$.

\par Since each agent needs to perform local particle smoothing before communication, it is necessary to consider the time complexity $\mathcal{O}(2TMn)$ and space complexity $\mathcal{O}(M+\frac{M}{{J_v}_{\text{max}}})$. It can be seen that the computational complexity of the proposed algorithm is greatly affected by the number of particles.

\subsection{Convergence Analysis}
\label{6.3}

\par The difficulty of the convergence analysis for Algorithm \ref{al4} mainly lies in analyzing the particle approximation $\bar{{\hat{Q}}}({\theta}, \theta_{k})$. In the subsequent discussion, the global consensus of the state estimation in Algorithm \ref{al3} will be analyzed. Furthermore, it will be illustrated that Algorithm \ref{al4} can identify the model parameters that satisfy the contraction stability constraint (\ref{eq5.2.7}).

\par To demonstrate that Algorithm \ref{al3} can achieve the global particle consensus among network agents, it is essential to redefine a vector $\mathbf{c}_{t,v}\in\mathbb{R}_{\lceil\frac{M}{{J_v}_{\text{max}}}\rceil\times V}$ for each agent $v$. This vector tracks the neighbors $\mathcal P_v$ that sent data to $v$ during the previous communication time. Initially, at time $0$, the $(v-1)\lceil\frac{M}{{J_v}_{\text{max}}}\rceil+1$ to $v\lceil\frac{M}{{J_v}_{\text{max}}}\rceil$ coordinates of this vector are set to 1, while the other coordinates are set to 0. Subsequently, at time $t$, it is updated as $\mathbf{c}_{t,v}=\sum_j \mathbf{c}_{t-1,j}$, $j\in\mathcal{P}_v$. Utilizing this vector, the particle transmission vector received by agent $v$ at time $t$ can be expressed as $ {Xc}_{t,v}=\sum_j \mathbf{c}_{t-1,j}{Xc}_{t-1,j}$, and the propagation number is defined as $n_{t,v}=\frac{\left\| \mathbf{c}_{t,v}\right\|_1}{\lceil\frac{M}{{J_v}_{\text{max}}}\rceil}$. When $\mathbf{c}_{t,v}$ approximates an all-1 vector $\mathbf{1}$, the average $\frac{Xc_{t,v}}{n_{t,v}}$ approaches the true global average. For a network consisting of $V$ agents, the relative error of agent $v$ at time $t$ is defined as $\sigma_{t,v}=\text{max}_k\left|\frac{{c}_{t,vk}}{\left\| \mathbf{c}_{t,v}\right\|_1}-\frac{1}{V} \right|=\left\|\frac{\mathbf{c}_{t,v}}{\left\| \mathbf{c}_{t,v}\right\|_1}-\frac{1}{V}\cdot \mathbf{1}\right\|_\infty$, where ${c}_{t,vk}$ is the $k$-th element of the vector $\mathbf{c}_{t,v}$.
\begin{thm}
	\label{thm3}
	The proposed PC-DPF algorithm (Algorithm \ref{al3}) can ensure that the relative error $e$ of the $\text{mean}(\{Xc_{t,v}\}_{v=1}^{V})$ and $\frac{Xc_{t,v}}{n_{t,v}}$ is less than $\delta$ with a probability of at least $1 - \frac{1}{V}$ at time $t=\log_2 V+\log_2 \frac{1}{\bar \delta}$.
\end{thm}
\begin{pf}
	Theorem 5.2 in \cite{FRIEZE198557} indicates that, for the considered message-passing scheme, the probability of information broadcasting to all agents within $(1 +\epsilon)(1 +(\gamma+ 1) \log 2) \log_2 n$ time steps is $1-n^{-\gamma}$. By choosing $\gamma=1$ and $\epsilon=\frac{1}{2}$, and defining $\tau=(2+3\log2)\log_2 2V$, in accordance with step 2.2 in Algorithm \ref{al3}, the propagation count in each broadcast is $\frac{1}{2}n_{t,v}$. Consequently, over $\tau$ transmissions, each agent's propagation count is at least $2^{-\tau}$, and this event occurs with a probability of at least $1-\frac{1}{2V}$.
	
	\par Define the function $\Phi_t=\sum_{v}\sum_{k}(c_{t,vk}-\frac{n_{t,v}}{V})^2$, where $\Phi_0\leqslant V$. It is easy to show that the conditional expectation of $\Phi_{t+1}$ at time $t+1$ satisfies $E[\Phi_{t+1}|\Phi_{t}]=(\frac{1}{2}-\frac{1}{2V})\Phi_{t}$ (refer to the proof of Lemma 2.3 in \cite{1238221}). Define the absolute error $\bar \delta=\delta^2\frac{1}{2V}2^{-2\tau}$, where $\delta$ is an arbitrarily small positive number. After $t=\log_2 V+\log_2 \frac{1}{\bar \delta}$ rounds of message passing, the expectation of $\Phi_{t}$ satisfies
	\begin{equation}
		\setlength{\abovedisplayskip}{3pt}
		\setlength{\belowdisplayskip}{3pt}
		\textstyle
		\begin{aligned}
			E[\Phi_{t}|\Phi_{0}\leqslant V]&\leqslant V\cdot2^{-t}\\
			&= V\cdot2^{-(\log_2 V+\log_2 \delta^{-2}+\log_2 2V+\log_2 2^{2\tau})}\\
			&=\bar \delta\notag
		\end{aligned}
	\end{equation}
 Substituting this into the Markov's inequality yields $p(\Phi_{t} \leqslant \delta^2 2^{-2\tau}) \geqslant 1 - \frac{1}{2V}$. Then, it can be derived that $p(|c_{t,vk}-\frac{n_{t,v}}{V}|\leqslant \delta 2^{-\tau}) \geqslant 1 - \frac{1}{2V}$.
 
 \par Applying Boole's inequality, it has $p(|\frac{ c_{t,vk}} {n_{t,v}}-\frac{1}{V}|\leqslant \delta ) \geqslant 1 - \frac{1}{V}$. Therefore, with a probability not less than $1 - \frac{1}{V}$, it can be asserted that, at time $t$, the relative error of agent $v$ satisfies $\sigma_{t,v}=\left\|\frac{\mathbf{c}_{t,v}}{\left\| \mathbf{c}_{t,v}\right\|_1}-\frac{1}{V}\cdot \mathbf{1}\right\|_\infty\leqslant \delta$. This completes the proof of the theorem.
 {\hfill $\square$}
\end{pf}

\par Based on the above theorem, the convergence of the global consensus state estimation particles can be guaranteed. Consequently, the identification of parameters by maximizing the likelihood function $L_{\theta}(X(T),Y(T))$ in the large-scale network comprising $V$ isomorphic agents is equivalent to identifying the parameter from the data of a single agent in Algorithm \ref{al4}. Since $Q({\theta}, \theta_{k})$ is replaced by the particle approximation $\bar{\hat{Q}}({\theta}, \theta_{k})$, the following lemma is needed.
\begin{lem}
	\label{l2}
	Consider the function $Q({\theta}, \theta_{k})$ defined by (\ref{eq3.1.9}) and its particle approximation $\bar{\hat{Q}}({\theta}, \theta_{k})$ defined by (\ref{eq5.1.4}) which is based on $M$ particles. Suppose that 
		\begin{equation}
			\label{eq7.1}
			\setlength{\abovedisplayskip}{3pt}
			\setlength{\belowdisplayskip}{3pt}
			\textstyle
			\begin{aligned}
				&p_{\theta} ({y}_v(t) | x_v(t) ) < \infty,\quad p_{\theta} (x_v(t+1) | x_v(t) ) < \infty,  \\
				&E\{|Q({\theta}, \theta_{k})|^4| Y_v(T)\} < \infty,
			\end{aligned}
		\end{equation}	
		hold for all ${\theta}$, $\theta_{k}\in\Theta$. Then with probability one
		\begin{equation}
			\label{eq7.2}
			\setlength{\abovedisplayskip}{3pt}
			\setlength{\belowdisplayskip}{3pt}
			\textstyle
			\underset{M\to\infty}{\lim} \bar{\hat{Q}}({\theta}, \theta_{k}) ={Q}({\theta}, \theta_{k}), 
			\forall {\theta}, \theta_{k}\in\Theta.
	\end{equation}
\end{lem}
\begin{pf}
	By application of Corollary 6.1 in \cite{4471882}.{\hfill $\square$}
\end{pf}

\par It can be shown that, for a sufficiently large number of particles $M$, the particle approximation $\bar{\hat{Q}}({\theta}, \theta_{k})$ for $Q({\theta}, \theta_{k})$ is both reasonable and accurate. Therefore, the combination of Lemma \ref{l2} with Theorem \ref{thm1} and Theorem \ref{thm5} can establish that $\hat{\theta}_{k+1}$ can be found such that $\bar{\hat{Q}}(\hat{\theta}_{k+1}, {\theta}_{k}) > \bar{\hat{Q}}({\theta}_{k},\theta_{k})$. That is, the iterative result $\hat{\theta}_{k+1}$ will be closer to the maximum likelihood estimate.

\par In conclusion, Theorem \ref{thm3} establishes the global consensus of state estimation particles for isomorphic network agents, enabling the distributed identification of the network parameter. The formulations presented in Theorem \ref{thm1} and Theorem \ref{thm5} lay a scientific foundation for the distributed parameter estimation using $\bar{\hat{Q}}({\theta}, \theta_{k})$. This section provides the theoretical guarantee for the proposed PC-DPEM algorithm to achieve distributed parameter identification of an isomorphic large-scale network.

\section{Simulation}
\label{Sec7}

\par In this section, the practicality and performance of Algorithm \ref{al4} are demonstrated through simulation examples. The simulation considers a complex large-scale network topology based on a synthetic network, generated using the Barab\'asi-Albert (BA) \cite{network} model with scale-free (SF) topology. The Barab\'asi-Albert diagram $(n, m)$ from networkX \cite{Hagberg2008ExploringNS} is utilized, where $n = 100$ and $m = 5$. Each link is set to be bidirectional, and then a portion of these unidirectional links is randomly deleted, ensuring that the network remains sparsely connected. The directed network has an average total degree of 5.1, a maximum degree of 40, and a minimum degree of 3. The maximum number of iterations in Algorithm \ref{al3} is set to 1000, and its iteration termination condition is determined by either the number of iterations or the convergence error. To test the effectiveness of Algorithm \ref{al4} in all the cases described below, the Monte Carlo study is performed using 100 different data realizations $Y_N$ of length $N = 100$. For each case, an estimate $\hat{\theta}$ is computed using Algorithm \ref{al4} with random initialization $\theta_{0}$ such that each entry of $\theta_{0}$ lies in an interval equal to $50\%$ of the corresponding entry in the true parameter vector $\theta^*$. In all cases, the proposed algorithm runs on an 11th-generation Intel Core i7-11700K processor of Windows 11. 

\subsection{General Nonlinear System}

\par Under the above network topology, the proposed PC-DPEM algorithm is used for the parameter identification of a network consisting of well-studied and challenging nonlinear subsystems\cite{SCHON201139} that interact with each other in the presence of uncertain interaction. The nonlinear and time-varying agent dynamics considered has the form of 
\begin{equation}
	\label{7.1}
	\setlength{\abovedisplayskip}{3pt}
	\setlength{\belowdisplayskip}{3pt}
	\textstyle
	\begin{aligned}
		&\textstyle{x}_v(t+1)=a{x}_v(t)+b\frac{x_v(t)}{1+x_v^2(t)}+c\cos(1.2t)\\
		&\textstyle\qquad\qquad\quad+\varepsilon_v(t)
		+\sum_{j\in\mathcal P_v}{g_j}(x_j(t)-x_v(t)),\\
		&\textstyle y_v(t)=dx_v^2(t)+\eta_{v}(t),\\
		&\textstyle\begin{bmatrix}
			\varepsilon_v(t)\\
			\eta_{v}(t)
		\end{bmatrix}\sim\mathcal{N}\left(\begin{bmatrix}
			0\\
			0
		\end{bmatrix},\begin{bmatrix}
			s & 0 \\
			0 & w
		\end{bmatrix}\right),
	\end{aligned}
\end{equation}
where the true parameters in this case are
\begin{equation}
	\label{7.2}
	\setlength{\abovedisplayskip}{3pt}
	\setlength{\belowdisplayskip}{3pt}
	\textstyle
	\begin{aligned}
		\theta^*&=\begin{bmatrix}a^* & b^* & c^* & d^* & s^* & w^*\end{bmatrix}\\
		&=\begin{bmatrix}0.5 & 25 & 8 & 0.05 & 0.5 & 1\end{bmatrix},
	\end{aligned}\notag
\end{equation}
and ${g_j}(x_j(t)-x_v(t))=\frac{1}{J_v}\sin(x_j(t)-x_v(t))$. 

\begin{figure}[htbp] 
	\centering
	\includegraphics[width=0.5\textwidth,trim=0 0 0 0,clip]{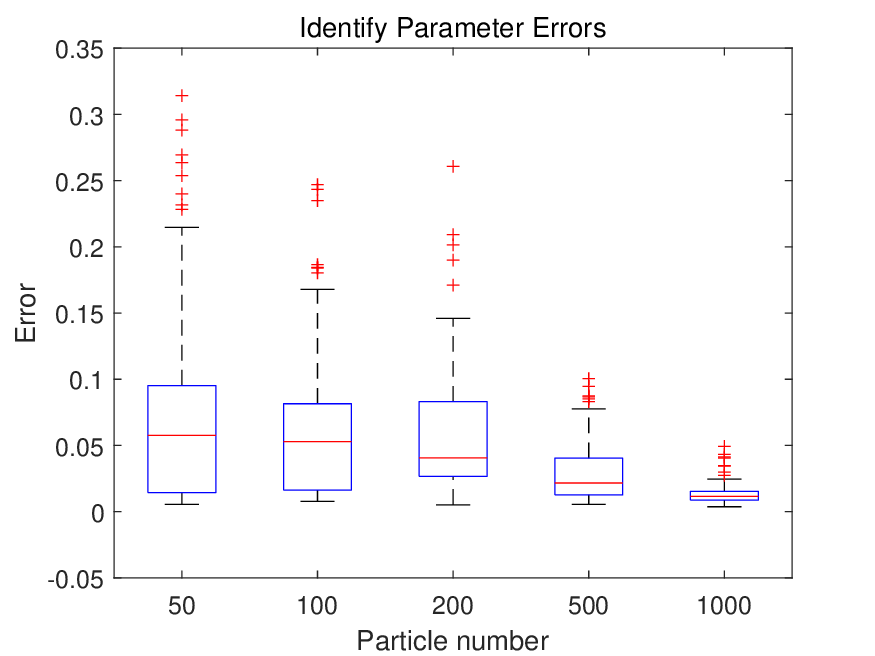}
	\caption{Parameter identification errors for different particle numbers.}
	\label{f2}
\end{figure}

\par The simulation example in this section takes into account the direct correlation between the number of particles and the identification accuracy and computational efficiency. Fig. \ref{f2} shows the parameter identification errors with 50, 100, 200, 500, and 1000 particles, respectively, where the error is calculated by $\frac{\|\hat\theta-\theta^*\|}{\|\theta^*\|}$. It can be seen that the parameter identification errors for different particle numbers are all less than 0.25, and with the increase in the number of particles, the parameter identification errors gradually decrease, leading to more accurate results. Table \ref{t1} lists the Monte Carlo test results with $1000$ particles, where the rightmost column presents the sample mean plus/minus sample standard deviation of the Monte Carlo trials. It can be observed that the parameter vector can be accurately estimated using a large number of particles.

\begin{table}[htbp]
	\caption{True and estimated values for isomorphic agent dynamics (\ref{7.1}).}
	\label{t1}
	\centering
	\def\temptablewidth{0.4\textwidth}
	\begin{tabular*}{\temptablewidth}{@{\extracolsep{\fill}}ccccc}
		\toprule[0.5pt]
		&Parameter&True&Estimated&\\
		\midrule[0.5pt]
		&a&0.5&$0.495\pm7.13\times10^{-4}$&\\
		&b&25&$24.9\pm0.119$&\\
		&c&8&$8.05\pm0.103$&\\
		&d&0.05&$0.053\pm1.45\times10^{-3}$&\\
		&s&0.5&$0.451\pm0.929$&\\
		&w&1&$1.18\pm0.432$&\\
		\bottomrule[0.5pt]
	\end{tabular*}
\end{table}

\par
\begin{figure}[htbp]   
	\centering
	\includegraphics[width=0.5\textwidth,trim=0 0 0 0,clip]{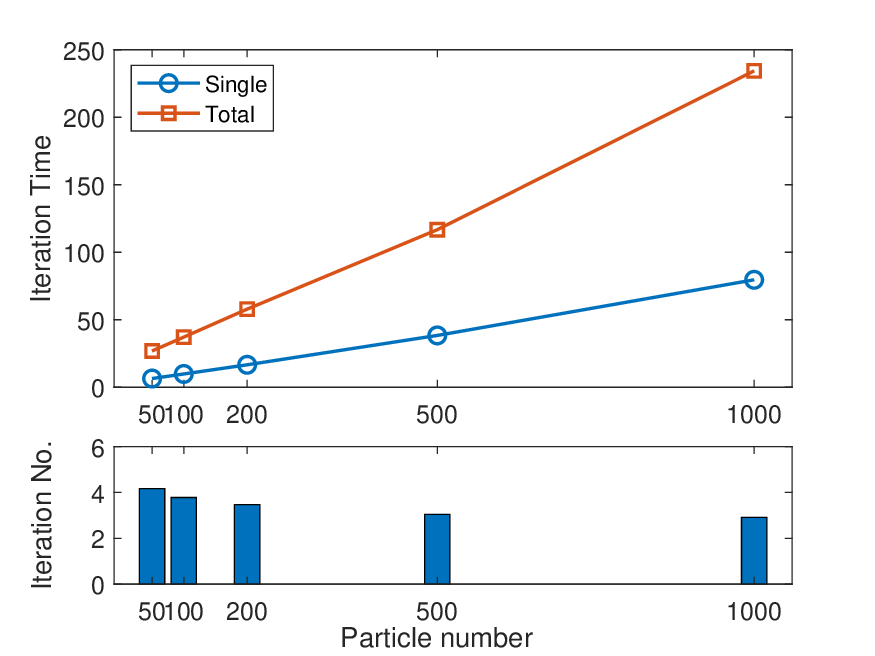}
	\caption{Iteration time and number of iterations for obtaining optimal parameters under different particle number conditions.}
	\label{f3}
\end{figure}

 As it is mentioned in Section \ref{Sec6} that the number of particles has a direct correlation with the communication overhead and computational complexity, this section further compares the single iteration time, total iteration time, and the number of iterations to run the proposed PC-DPEM algorithm with different particle numbers in Fig. \ref{f3}. From the figure, it is evident that as the number of particles increases, the single iteration time grows from $4.6s$ to $54.9s$, and the total iteration time increases from $19.5s$ to $148.7s$. This is attributed to the increased involvement of particle states in transmission as the number of particles increases. Similarly, the number of iterations for different agents to communicate with neighboring agents to reach a consensus also increases. However, the analysis in Section \ref{6.3} indicates that larger particle numbers bring the identification result much closer to the real value. This is reflected in the decrease in the number of iterations, dropping from $4.37$ for 50 particles to $2.71$ for $1000$ particles.

\par
\begin{figure}[htbp]
	\centering
	\begin{subfigure}[b]{0.236\textwidth}
		\includegraphics[width=\textwidth,trim=5 0 38 8,clip]{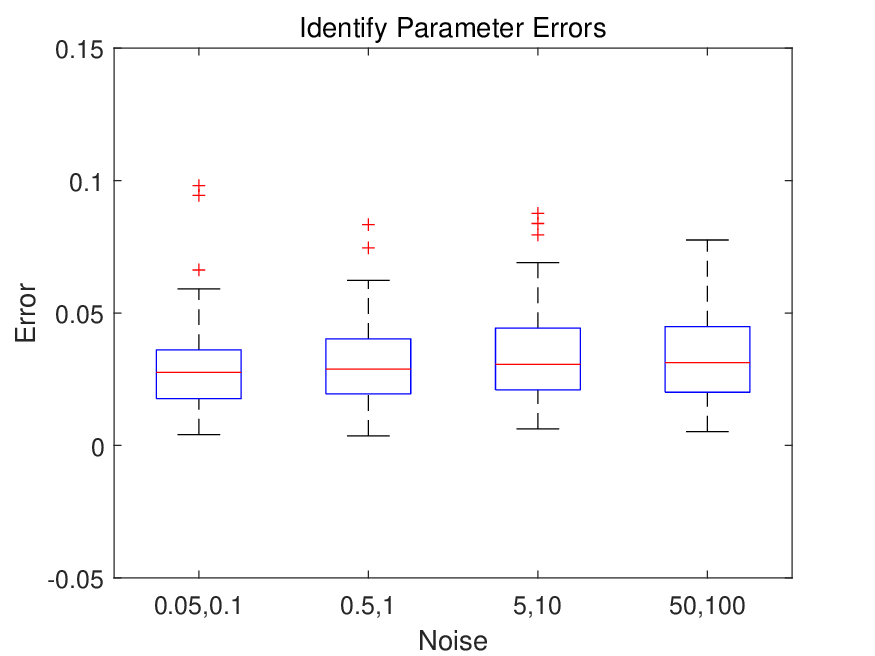}
		\vspace{-0.5cm}
		\caption{ }
		\label{f4.1}
	\end{subfigure}
	\begin{subfigure}[b]{0.236\textwidth}
		\includegraphics[width=\textwidth,trim=15 0 28 8,clip]{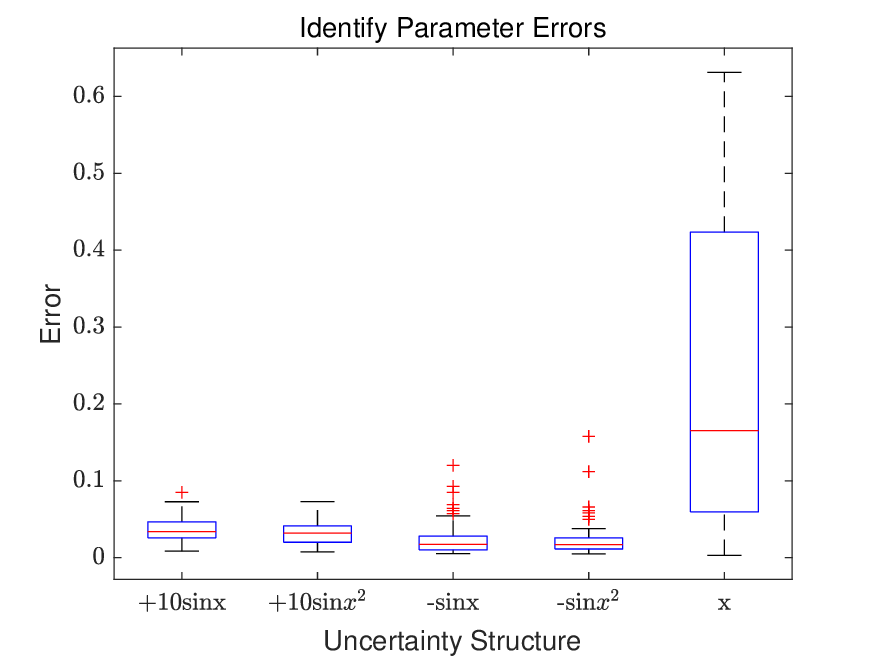}
		\vspace{-0.5cm}
		\caption{ }
		\label{f4.2}
	\end{subfigure}
	\caption{Parameter identification errors for different conditions. (a) Four different noise conditions. (b) Different structures of unknown state interactions.}
\end{figure}

 In addition, this section validates the parameter identification results under different noise conditions and various network uncertain interactions $\sum_{j\in\mathcal P_v}{g_j}(x_j(t),x_v(t))$ based on the 1000-particle condition. Fig. \ref{f4.1} compares the parameter identification results in the presence of noise conditions with $\mathcal{N}(0.05,0.1)$, $\mathcal{N}(0.5,1)$, $\mathcal{N}(5,10)$, and $\mathcal{N}(50,100)$, respectively. It can be observed that, although the parameter identification error slightly increases as the noise conditions become worsen, the overall error remains below $0.1$. This demonstrates that the proposed algorithm can achieve accurate parameter identification even under harsh noise conditions. 

\par Fig. \ref{f4.2} compares the model parameter identification results in the presence of unknown state interaction between agents shown in Table \ref{t2}. It can be observed that the proposed algorithm yields identification results with errors less than $0.6$ even though the interactions among the agents are unknown. This section also verifies the parameter identification results in the presence of the above uncertain interactions and without a priori stability constraints. Note that the iterative process diverges directly under these conditions and is therefore not shown in Fig. \ref{f4.2}.

	\begin{table}[htbp]
	\caption{Uncertain agent state interaction structure ${g_j}(x_j(t)-x_v(t))$.}
	\label{t2}
	\centering
	\def\temptablewidth{0.4\textwidth}
	\begin{tabular*}{\temptablewidth}{@{\extracolsep{\fill}}ccc}
		\toprule[0.5pt]
		&${g_j}(x_j(t),x_v(t))$ functions&\\
		\midrule[0.5pt]
		&${10}/{J_v}\sin(x_j(t)-x_v(t))$&\\
		&${10}/{J_v}\sin(x_j(t)-x_v(t))^2$&\\
		&$-{1}/{J_v}\sin(x_j(t)-x_v(t))$&\\
		&$-{1}/{J_v}\sin(x_j(t)-x_v(t))^2$&\\
		&${1}/{J_v}(x_j(t)-x_v(t))$&\\
		\bottomrule[0.5pt]
	\end{tabular*}
\end{table}
	
\subsection{Complex Network Dynamics Related to Social and Brain Systems}
	
\par In order to demonstrate the effectiveness of the method proposed in this study for the identification of complex network dynamics, the proposed Algorithm \ref{al4} is applied to identify isomorphic complex networks related to social and brain systems. The topological relationships of the agents in the network still follow the BA model with the SF topology used in the previous subsection. The following simulations are all based on the 500-particle condition. The conditions of random initialization change for the parameters in the network, while the length of the identification data remains unchanged from the settings in the previous subsection.

\subsubsection{Gene regulation dynamics}

\par This subsection first demonstrates the effectiveness of the proposed algorithm in inferring gene regulation (GR) dynamics. Gene activity data can be generated according to the equation in \cite{20038296}. In order to highlight the applicability of the proposed algorithm to conditions such as incomplete state observability and the presence of noise, the GR dynamics has the following form
\begin{equation}
	\label{6.2.1}
	\setlength{\abovedisplayskip}{3pt}
	\setlength{\belowdisplayskip}{3pt}
	\textstyle
	\begin{aligned}
		\textstyle\frac{d x_v}{dt}&\textstyle=ax_v+\varepsilon_v+\epsilon_{vj}\sum_{j\in\mathcal P_v}\frac{(x_j-x_v)^\alpha}{(x_j-x_v)^\alpha+1}\\
		y&\textstyle=bx_v+\eta_v,
	\end{aligned}
\end{equation}
where $x_v(t)$ is the concentration of gene $v$ at time $t$, $\alpha = 2$ denotes Hill's coefficient, the true parameters in this dynamics are
\begin{equation}
	\setlength{\abovedisplayskip}{3pt}
	\setlength{\belowdisplayskip}{3pt}
	\textstyle
		\theta^*_v=\begin{bmatrix}a^* & b^* & s^* & w^*\end{bmatrix}
		=\begin{bmatrix}-0.2 & 1 & 0.001 & 0.01\end{bmatrix}.\notag
\end{equation} 
The $\frac{(x_j-x_v)^\alpha}{(x_j-x_v)^\alpha+1}$ is the interaction function between gene $v$ and its neighbor $\mathcal P_v$. The parameter $\epsilon_{vj}$ denotes the strength of the gene $j$'s regulation of gene $v$, which is set to $0.05$ in the network. For verification purposes, the initial value of the state $x_{v0}$ in the dynamics (\ref{6.2.1}) of each of the 100 agents is set to $x_{v0}=1$.

	\begin{table}[htbp]
	\caption{True, estimated values of isomorphic gene generation dynamics (\ref{6.2.1}), and the comparison with the identification results of \cite{Gao}.}
	\label{t3}
	\centering
	\def\temptablewidth{0.46\textwidth}
	\setlength{\tabcolsep}{-4pt}
	\begin{tabular*}{\temptablewidth}{@{\extracolsep{\fill}}cccccc}
		\toprule[0.5pt]
		&Parameter&True&Estimated&Comparison&\\
		\midrule[0.5pt]
		&a&-0.2&-$0.2003\!\pm\!1.23\!\times\!10^{-4}$&-0.197&\\
		&b&1&$1.0060\!\pm\!1.26\!\times\!10^{-4}$&-&\\
		&s&0.001&$0.00179\!\pm\!0.64$&-&\\
		&w&0.01&$0.00998\!\pm\!0.59$&-&\\
		\bottomrule[0.5pt]
	\end{tabular*}
	\end{table}

\par
\begin{figure}[htbp]
		\centering
		\begin{subfigure}[b]{0.236\textwidth}
			\includegraphics[width=\textwidth,trim=20 5 28 5,clip]{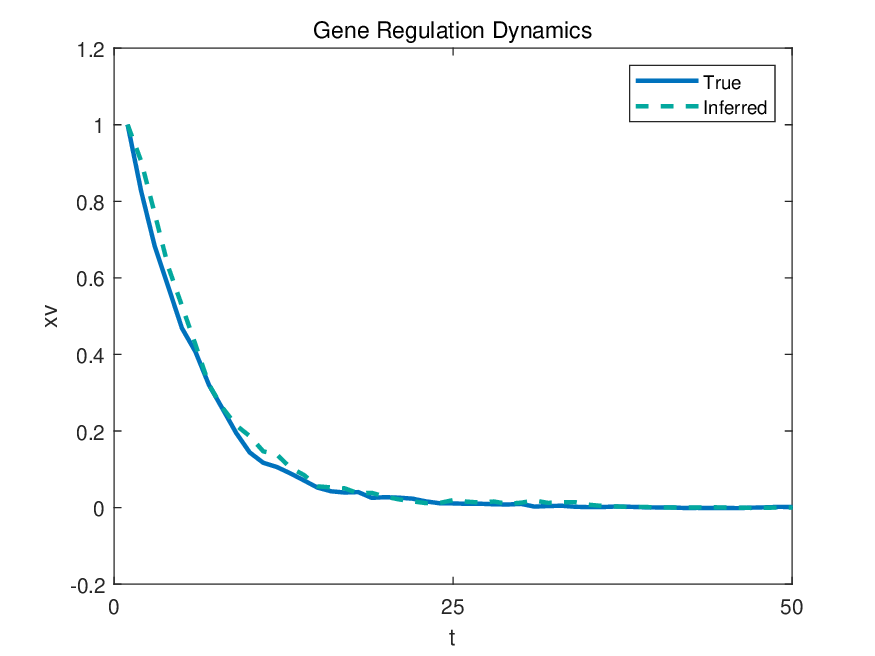}
			\vspace{-0.5cm}
			\caption{ }
			\label{f5.1}
		\end{subfigure}
		\begin{subfigure}[b]{0.236\textwidth}
			\includegraphics[width=\textwidth,trim=10 3 28 5,clip]{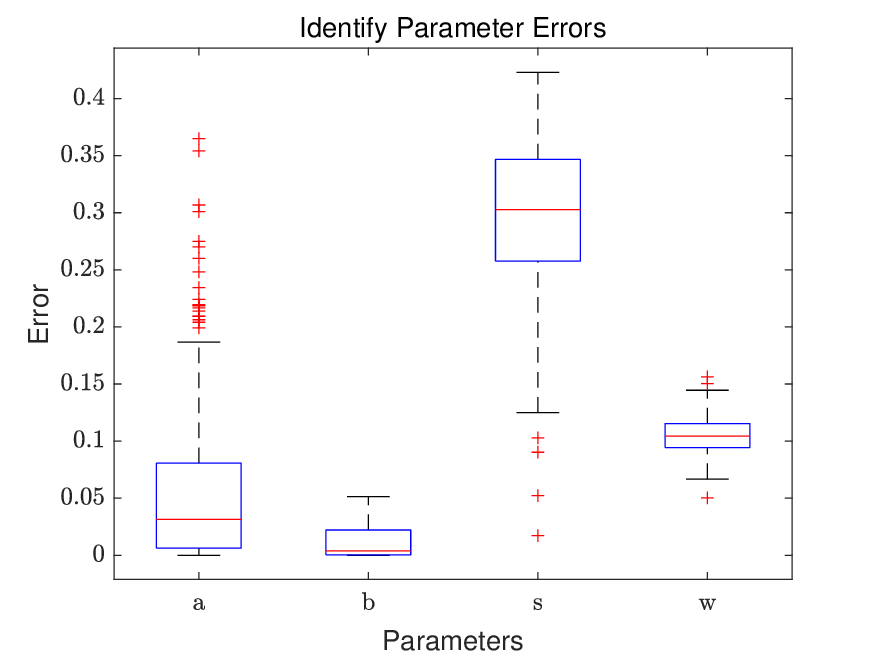}
			\vspace{-0.5cm}
			\caption{ }
			\label{f5.2}
		\end{subfigure}
		\caption{(a) True and inferred gene activity trajectories for gene $v$. (b) Parameter identification errors for dynamics (\ref{6.2.1}).}
	\end{figure}

 Table \ref{t3} shows the identification results of Algorithm \ref{al4} for unknown parameters of dynamic (\ref{6.2.1}). The rightmost parameter is the identification result obtained by \cite{Gao} without considering noise interference, and the system state and the interaction relationship between different genes are fully known. As can be seen from the data in the table, the results of parameter identification obtained by the Algorithm \ref{al4} under complex uncertainty conditions are similar to those obtained by \cite{Gao} under more ideal conditions. The results in the table are considered to be successful since the identification error for each parameter does not exceed $0.4$ (as shown in Fig. \ref{f5.2}) with a $50\%$ range of random initialization and a small number of outliers. The activity trajectories of gene $v$ generated from the true and inferred equations are shown in Fig. \ref{f5.1}.

\subsubsection{FitzHugh-Nagumo dynamics}

\par This subsection also tests the algorithm by applying it to neuronal activity data generated by FitzHugh-Nagumo (FHN) dynamics \cite{FITZHUGH1961445} on a BA network. Taking into account the presence of noise and the partial observability of the state, the model of the FHN neural network has the following form  
\begin{equation}
	\label{6.2.3}
	\setlength{\abovedisplayskip}{3pt}
	\setlength{\belowdisplayskip}{3pt}
	\textstyle
	\begin{aligned}
		\textstyle\frac{d x_{v1}}{dt}&\textstyle=ax_{v1}+bx_{v1}^3+cx_{v2}+\varepsilon_v\\
		&\textstyle\qquad\quad+\epsilon_{vj}\sum_{j\in\mathcal P_v}\frac{x_{j1}-x_{v1}}{J_v}
		\\\textstyle\frac{d x_{v2}}{dt}&\textstyle=d+ex_{v1}+fx_{v2}+\varepsilon_v\\
		y&\textstyle=x_{v1}+x_{v2}+\eta_v.
	\end{aligned}
\end{equation}
FHN dynamics capture the firing behaviors of neurons, which have two components. The first component $x_{v1}$ represents the membrane potential containing the interaction dynamics of the self $v$ and the neighboring neuron $\mathcal P_v$, where $J_v$ is the number of $j$, and $\epsilon_{vj}= 1$. The second component $x_{v2}$ represents the recovery variables. The true parameters in this dynamics are
\begin{equation}
	\setlength{\abovedisplayskip}{3pt}
	\setlength{\belowdisplayskip}{3pt}
	\begin{aligned}
		\textstyle
		\theta^*_v&\textstyle=\begin{bmatrix}a^* & b^* & c^* & d^* & e^* & f^* & s^* & w^*\end{bmatrix}\\
		&\textstyle=\begin{bmatrix}1 &-1 &-1 &0.28 & 0.5 & -0.04 & 0.05& 0.1\end{bmatrix}.\notag
	\end{aligned}
\end{equation} 
For verification purposes, the initial value of the state $x_{v0}$ in the dynamics (\ref{6.2.1}) of each of the 100 agents is set to $x_{v0}=[0.87609,-3.5091]$.

	\begin{table}[htbp]
		\setlength{\tabcolsep}{-1pt}
		\caption{True, estimated values of isomorphic FHN dynamics (\ref{6.2.3}), and the comparison with the identification results of \cite{Gao}.}
		\label{t4}
		\centering
		\def\temptablewidth{0.46\textwidth}
		\begin{tabular*}{\temptablewidth}{@{\extracolsep{\fill}}cccccc}
			\toprule[0.5pt]
			&Parameter&True&Estimated&Comparison&\\
			\midrule[0.5pt]
			&a&1&$0.9931\!\pm\!0.011$&0.989&\\
			&b&-1&-$1.0171\!\pm\!0.016$&-0.993&\\
			&c&-1&-$1.0086\pm0.003$&-0.996&\\
			&d&0.28&$0.2814\pm0.006$&0.279&\\
			&e&0.5&$0.4998\pm0.213$&0.499&\\
			&f&-0.04&-$0.0402\pm0.113$&-0.040&\\
			&s&0.05&$0.0466\!\pm\!0.629$&-&\\
			&w&0.1&$0.1060\!\pm\!0.514$&-&\\
			\bottomrule[0.5pt]
		\end{tabular*}
	\end{table}

\par
\begin{figure}[htbp]
		\centering
		\begin{subfigure}[b]{0.236\textwidth}
			\includegraphics[width=\textwidth,trim=10 0 38 6,clip]{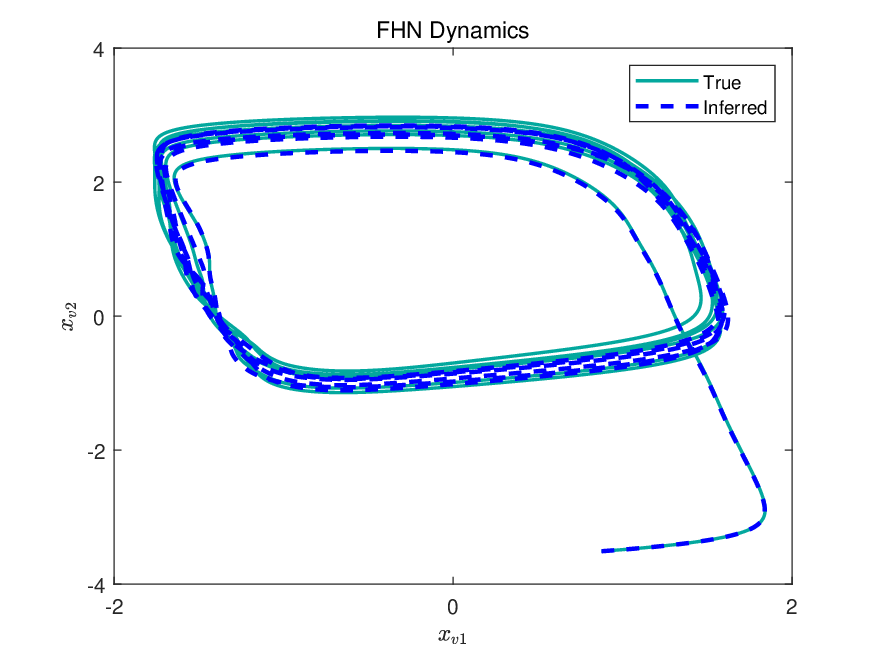}
			\vspace{-0.5cm}
			\caption{ }
			\label{f6.1}
		\end{subfigure}
		\begin{subfigure}[b]{0.236\textwidth}
			\includegraphics[width=\textwidth,trim=10 0 38 6,clip]{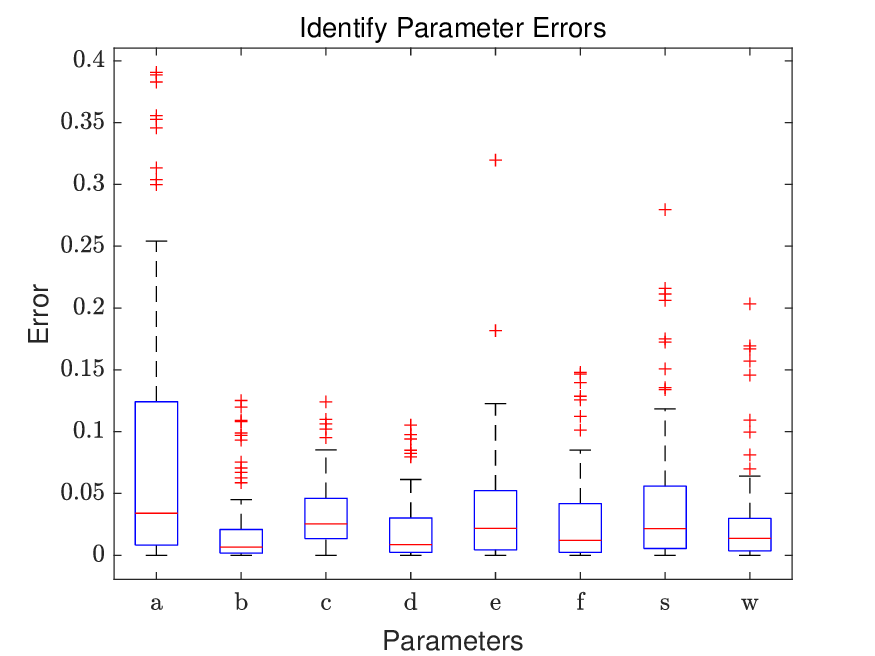}
			\vspace{-0.5cm}
			\caption{ }
			\label{f6.2}
		\end{subfigure}
		\caption{(a) True and inferred gene activity trajectories of neurons $v$. (b) Parameter identification errors for dynamics (\ref{6.2.3}).}
	\end{figure}

 Table {\ref{t4}} shows the results of Algorithm \ref{al4} for the identification of the unknown parameters of the dynamics (\ref{6.2.3}). The rightmost parameter is still the identification result obtained by \cite{Gao} without considering the noise interference and the system state and the interaction dynamics between different neurons are completely known. As can be seen from the data, the parameter identification results obtained by Algorithm \ref{al4} under complex uncertainty conditions are similar to those obtained by \cite{Gao} under more ideal conditions, and the absolute parameter identification error 0.0206 obtained by Algorithm \ref{al4} is smaller than that of 0.03 by \cite{Gao}. The results in the table are considered to be successful because the identification error for each parameter does not exceed $0.25$ under a random initialization range of $50\%$ and a small number of outliers (as shown in Fig. \ref{f6.2}). The activity trajectories of neurons $v$ generated from the true and inferred equations are shown in Fig. \ref{f6.1}.
	
\section{Conclusion}
	\label{Sec8}

\par In this paper, a particle consensus-based distributed particle EM algorithm has been presented by integrating DPF with EM algorithm under a priori contraction stability constraints. The algorithm accomplishes parameter identification for isomorphic large-scale nonlinear network dynamics, even under conditions of incomplete state observability, challenging noise environments, and the existence of unknown interactions among network agents. Performance analysis of the proposed method and validation through simulation results confirm its effectiveness in identifying parameters for stable nonlinear isomorphic networks. The next step is planned to investigate the identification of heterogeneous large-scale networks, as well as to improve the computational efficiency of the current algorithm.


\bibliographystyle{unsrt}        
\bibliography{autosam}           



\end{document}